\documentclass[11pt]{article}

\usepackage[utf8]{inputenc}
\usepackage[T1]{fontenc}
\usepackage{amsmath}
\usepackage{amssymb}
\usepackage{amsthm}
\usepackage[margin=1in]{geometry}
\usepackage{enumitem}
\usepackage{booktabs}
\usepackage{array}
\usepackage{graphicx}
\usepackage{xcolor}
\usepackage{hyperref}
\usepackage[protrusion=true,expansion=false]{microtype}
\hypersetup{
  colorlinks=true,
  linkcolor=blue!50!black,
  citecolor=blue!50!black,
  urlcolor=blue!50!black,
  pdftitle={Decimal Expansion of lambda\_infinity},
  pdfauthor={Alessandro Munari}
}

\title{$\lambda_\infty$: A New Mathematical Constant from the Spectral Theory of the Prime LCM Matrix}
\author{Alessandro Munari\thanks{Email: \href{mailto:almunari@icloud.com}{almunari@icloud.com}}}
\date{June 2026}

\begin{document}

\maketitle

\begin{abstract}
We introduce a new mathematical constant
\[\lambda_\infty = 0.674036183193696139936660007576508455780\ldots \quad \text{(OEIS A396695)},\]
defined as the unique solution in $\bigl(\tfrac{1}{4},+\infty\bigr)$ (Theorem 3.1) of
\[h(x) := \sum_{p,\text{prime}} \frac{1}{xp^2-p+1} = 1.\]
This equation arises as the $n\to\infty$ limit of the secular equation of rank-$n$ truncations of the infinite prime LCM matrix $\mathcal{L}[p_i,p_j]=1/\operatorname{lcm}(p_i,p_j)$ (Lemma 5.3), where $\mathcal{L}=D+vv^T$ for a diagonal part $D$ and $v=(1/p)_{p\in\mathcal{P}}$. Viewed as a compact self-adjoint operator on $\ell^2(\mathcal{P})$, $\mathcal{L}$ has spectral radius $\rho(\mathcal{L})=\lambda_\infty$ (Theorem 5.4), distinguishing it from the spectral radius $P(2)\approx 0.4522$ of the underlying rank-one operator $vv^T$ (equivalently, the spectral limit $\lim_{N\to\infty}\lambda_{\max}(W_N)/N$ of the integer prime divisor matrices; see Remark~5.1). The corresponding integer LCM matrix satisfies $\lambda_{\max}(W_N^*)/N\to\lambda_\infty$ (Theorem 5.8), a prime-indexed counterpart to the classical identity $\lambda_{\max}(M_N)/N\to\zeta(2)=\pi^2/6$ for the integer divisor matrix.

We prove that $h$ is real-analytic, strictly decreasing and strictly convex on $\bigl(\tfrac{1}{4},+\infty\bigr)$, ensuring existence and uniqueness of $\lambda_\infty$. The exact power-series expansion $h(x)=\sum_{k\ge 0}Q_k/x^{k+1}$, with coefficients expressible via the prime zeta function $P(s)=\sum_p p^{-s}$, yields an efficient Newton--Raphson algorithm with geometric prime-zeta tail correction for high-precision computation.

We compute 500 decimal digits of $\lambda_\infty$, certified by rigorous error bounds (Propositions~7.1--7.2) and independently verified by six computational runs, including a fully independent recomputation (Run~E, PARI-GP; ancillary file \texttt{run\_E\_pari\_recompute.gp}) and a machine-verifiable Arb interval certificate (Run~F; Proposition~7.2) covering 505 digits. Extensive PSLQ and LLL searches find no minimal polynomial of degree $\le 8$ satisfied by $\lambda_\infty$ (at 560 decimal digits), and no integer relation against catalogs of up to 31 classical constants (at 160--300 decimal digits). What, if anything, $\lambda_\infty$ is algebraically, remains open.
\end{abstract}

\tableofcontents

\section{Introduction}

The prime zeta function $P(s)=\sum_{p,\text{prime}}p^{-s}$ and its special values sit at the heart of analytic number theory. The value $P(2)=\sum_p p^{-2}\approx 0.4522$ (OEIS A085548~\cite{ref3}) is the prime-indexed counterpart of Euler's $\pi^2/6=\sum_{n\ge 1}n^{-2}$ (OEIS A013661~\cite{ref4}); both arise as the solutions of the implicit equations $\sum_p 1/(xp^2)=1$ and $\sum_{n\ge 1}1/(xn^2)=1$ respectively.

In this paper we study a deformed version of the prime equation: we include a linear correction $-p+1$ in each denominator, yielding
\[h(x) := \sum_{p,\text{prime}} \frac{1}{xp^2-p+1} = 1. \tag{1}\]
For every $x>\tfrac{1}{4}$ each denominator satisfies $xp^2-p+1>0$, since $(p-1)/p^2\le\tfrac{1}{4}$ for all primes (with equality at $p=2$), and the series $h(x)$ converges absolutely: for primes $p>2/x$ one has $xp^2-p+1>xp^2/2$, so those terms are dominated by $2/(xp^2)$, giving convergence by comparison with $(2/x)P(2)<\infty$; uniform convergence on compact subsets of $\bigl(\tfrac{1}{4},+\infty\bigr)$ follows similarly (Proposition 2.2). We prove that equation (1) has a unique solution in $\bigl(\tfrac{1}{4},+\infty\bigr)$, which we call $\lambda_\infty$.

The perturbation $-p+1$ is not negligible in its cumulative effect: it shifts the root from $P(2)\approx 0.4522$ to $\lambda_\infty\approx 0.6740$, a displacement of $\lambda_\infty - P(2)\approx 0.2218\approx 0.49\,P(2)$. The shift is almost entirely a low-prime effect. For each individual large prime the term-wise correction is asymptotically small, with $(-p+1)/(xp^2)\sim -1/(xp)\to 0$ as $p\to\infty$; but for $p=2$ the denominator decreases from $4\lambda_\infty\approx 2.696$ to $4\lambda_\infty-1\approx 1.696$, a relative reduction of $1/(4\lambda_\infty)\approx 0.371$ in the denominator of the dominant ($p=2$) term.

The constant $\lambda_\infty$ arose in the author's investigation of the spectral theory of prime divisor matrices. Let $\mathcal{P}=\{2,3,5,7,\ldots\}$ and let $\ell^2(\mathcal{P})$ be the Hilbert space of square-summable sequences indexed by $\mathcal{P}$. The matrix $\mathcal{L}[p_i,p_j]=1/\operatorname{lcm}(p_i,p_j)$ defines a bounded self-adjoint operator on $\ell^2(\mathcal{P})$. It decomposes as $\mathcal{L}=D+vv^T$, where $v=(p^{-1})_{p\in\mathcal{P}}\in\ell^2(\mathcal{P})$ ($\|v\|^2=P(2)<\infty$) and $D=\operatorname{diag}\bigl((p-1)/p^2\bigr)_{p\in\mathcal{P}}$ is the diagonal operator with $\|D\|_{\mathrm{op}}=\tfrac{1}{4}$. A direct computation yields
\[\|\mathcal{L}\|_{\mathrm{HS}}^2 = P(2) + P(2)^2 - P(4) < \infty\]
(see Remark 5.7(i)), so $\mathcal{L}$ is Hilbert--Schmidt, hence compact; being also self-adjoint, it has essential spectrum $\sigma_{\mathrm{ess}}(\mathcal{L})=\{0\}$. Its nonzero spectrum consists of countably many positive eigenvalues accumulating at $0$; by the rank-one perturbation interlacing theorem applied to the decomposition $\mathcal{L}=D+vv^T$ (the secular equation $(\star)$ of Lemma 5.3 shows that each truncation $\mathcal{L}_n=D_n+v_nv_n^T$ has exactly one eigenvalue exceeding $\|D_n\|_{\mathrm{op}}\le\tfrac{1}{4}$, namely $\rho(\mathcal{L}_n)$, while all remaining $n-1$ eigenvalues lie in $(0,\tfrac{1}{4}]$), all eigenvalues of $\mathcal{L}$ lie in $(0,\tfrac{1}{4}]$ except $\lambda_\infty$, which is the unique eigenvalue exceeding $\|D\|_{\mathrm{op}}=\tfrac{1}{4}$ and coincides with the spectral radius $\rho(\mathcal{L})$. Theorem 5.4 establishes this identification rigorously. While $\lambda_\infty$ is distinct from $\lim_{N\to\infty}\rho(W_N/N)$ --- where $W_N[i,j]=\lfloor N/(p_ip_j)\rfloor$ and the limit equals $P(2)$ (see Remark 5.1) --- it has a precise spectral identity. The term $-p+1$ in equation (1) is the exact algebraic consequence of $\operatorname{lcm}(p,p)=p$ (rather than $p^2$), as Theorem 5.4 establishes.

\textit{Notation.} Throughout, $p$ always denotes a prime number, and sums over $p$ always range over all primes $2,3,5,7,\ldots$ unless stated otherwise. We write $\mathcal{P}=\{2,3,5,7,\ldots\}$ for the set of all primes, $P(s)=\sum_p p^{-s}$ for the prime zeta function, and $\zeta(s)=\sum_{n\ge 1}n^{-s}$ for the Riemann zeta function.

\section{The function $h$: domain and analytic properties}

\textbf{Definition 2.1.} For $x\in\bigl(\tfrac{1}{4},+\infty\bigr)$ define
\[h(x)=\sum_{p,\text{prime}}\frac{1}{xp^2-p+1}.\]

\textbf{Proposition 2.2.} \textit{The function $h$ is well-defined, real-analytic, and positive on $\bigl(\tfrac{1}{4},+\infty\bigr)$.}

\textit{Proof.} The denominator $xp^2-p+1$ vanishes at $x=(p-1)/p^2$. For every prime $p$, the inequality
\[\frac{p-1}{p^2}\le\frac{1}{4}\]
is equivalent to $4(p-1)\le p^2$, i.e. $(p-2)^2\ge 0$, which holds for all real $p$ with equality precisely at $p=2$. Hence
\[\sup_{p,\,\text{prime}}\frac{p-1}{p^2}=\frac{2-1}{4}=\frac{1}{4}.\]
Thus for every prime $p$ and every $x>\tfrac{1}{4}$ we have $xp^2-p+1>0$, so every term is well-defined and positive.

For absolute and uniform convergence on any compact subinterval $[a,b]\subset\bigl(\tfrac{1}{4},+\infty\bigr)$, let $p_0=\lceil 2/a\rceil$. For every prime $p\ge p_0$ and every $x\in[a,b]$, one has $xp^2-p+1\ge ap^2-p+1\ge\tfrac{a}{2}p^2$ (since $p\ge p_0$ implies $\tfrac{a}{2}p\ge 1$), so
\[\frac{1}{xp^2-p+1}\le\frac{2}{ap^2},\]
a majorant uniform in $x\in[a,b]$. Since $\tfrac{2}{a}\sum_{p}p^{-2}=\tfrac{2}{a}P(2)<\infty$, the Weierstrass M-test gives absolute and uniform convergence on $[a,b]$.

For real-analyticity, fix any $x_0\in\bigl(\tfrac{1}{4},+\infty\bigr)$ and let $a=\tfrac{1}{2}(x_0+\tfrac{1}{4})>\tfrac{1}{4}$. For complex $z$ with $\operatorname{Re}(z)\ge a$, each denominator satisfies
\[|zp^2-p+1|\ge\operatorname{Re}(z)p^2-p+1\ge ap^2-p+1>\tfrac{a}{2}p^2\]
for $p\ge p_0=\lceil 2/a\rceil$. For the finitely many primes $p<p_0$, the denominator $zp^2-p+1$ has its unique zero at $z=(p-1)/p^2\le\tfrac{1}{4}$, which lies outside $\{\operatorname{Re}(z)>\tfrac{1}{4}\}$; on any compact subset $K\subset\{z\in\mathbb{C}:\operatorname{Re}(z)>\tfrac{1}{4}\}$ each such denominator is therefore bounded away from zero, so the corresponding finitely many terms are holomorphic and uniformly bounded on $K$, raising no convergence issue. Hence the series converges absolutely and uniformly on any compact subset of $\{z\in\mathbb{C}:\operatorname{Re}(z)>\tfrac{1}{4}\}$. In particular, the closed disc of radius $r=\tfrac{1}{2}(x_0-\tfrac{1}{4})$ centred at $x_0$ is contained in $\{\operatorname{Re}(z)>\tfrac{1}{4}\}$, and the uniform limit of holomorphic functions is holomorphic by the Weierstrass theorem, establishing real-analyticity at $x_0$. $\square$

\textbf{Proposition 2.3.} \textit{The function $h$ is strictly decreasing and strictly convex on $\bigl(\tfrac{1}{4},+\infty\bigr)$, with derivatives}
\[h'(x)=-\sum_p\frac{p^2}{(xp^2-p+1)^2}, \tag{2}\]
\[h''(x)=2\sum_p\frac{p^4}{(xp^2-p+1)^3}. \tag{3}\]
\textit{Both series converge absolutely and uniformly on compact subsets of $\bigl(\tfrac{1}{4},+\infty\bigr)$.}

\textit{Proof.} Each term $g_p(x)=(xp^2-p+1)^{-1}$ is real-analytic with $g_p'(x)=-p^2(xp^2-p+1)^{-2}<0$ and $g_p''(x)=2p^4(xp^2-p+1)^{-3}>0$. On any compact $[a,b]\subset\bigl(\tfrac{1}{4},+\infty\bigr)$, set $p_0=\lceil 2/a\rceil$. For every prime $p\ge p_0$ one has $ap^2-p\ge\tfrac{a}{2}p^2$, giving
\[\sum_{p\ge p_0}|g_p'(x)|\le\frac{4}{a^2}\sum_{p\ge p_0}p^{-2}<\infty.\]
The finitely many terms $g_p'$ (and $g_p''$) with prime $p<p_0$ are each $C^\infty$ on $[a,b]$ and bounded on the compact interval; their finite sum contributes a bounded differentiable function and raises no convergence issue.

An identical argument applies to the series for $h''$. Uniform convergence of the derived series justifies termwise differentiation, giving (2) and (3). Since every summand in (2) is negative (resp. in (3) is positive), we have $h'<0$ and $h''>0$ pointwise. $\square$

\textbf{Proposition 2.4.} \textit{The boundary behaviour of $h$ is:}

\textit{(i) $h(x)\to+\infty$ as $x\to\bigl(\tfrac{1}{4}\bigr)^+$.}

\textit{(ii) $h(x)\to 0$ as $x\to+\infty$.}

\textit{Proof.} (i) The $p=2$ term is $1/(4x-1)\to+\infty$ as $x\to\bigl(\tfrac{1}{4}\bigr)^+$. It remains to show that the residual sum $R(x):=\sum_{p\ge 3}(xp^2-p+1)^{-1}$ is bounded as $x\to\bigl(\tfrac{1}{4}\bigr)^+$. For every prime $p\ge 3$ and every $x\in\bigl(\tfrac{1}{4},\tfrac{1}{2}\bigr]$,
\[xp^2-p+1\;\ge\;\tfrac{1}{4}p^2-p+1\;=\;\frac{p^2-4p+4}{4}\;=\;\frac{(p-2)^2}{4},\]
so each term satisfies $1/(xp^2-p+1)\le 4/(p-2)^2$. The substitution $n:=p-2$ maps $\{p\ge 3:\text{$p$ prime}\}$ injectively into $\mathbb{N}_+$, so each term $4/(p-2)^2=4/n^2$ is a distinct summand of the convergent series $4\sum_{n\ge 1}n^{-2}=4\cdot\pi^2/6<\infty$; hence $\sum_{p\ge 3}4/(p-2)^2\le 4\sum_{n\ge 1}n^{-2}<\infty$. By the dominated convergence theorem (counting measure on the primes, with integrable majorant $4/(p-2)^2$), $R(x)$ converges to the finite limit $\sum_{p\ge 3}4/(p-2)^2$ as $x\to\bigl(\tfrac{1}{4}\bigr)^+$, and in particular remains bounded on $\bigl(\tfrac{1}{4},\tfrac{1}{2}\bigr]$. Therefore $h(x)=1/(4x-1)+R(x)\to+\infty$.

(ii) For fixed $p$ the term $(xp^2-p+1)^{-1}\to 0$ as $x\to+\infty$. For $x\ge 1$, each term satisfies
\[\frac{1}{xp^2-p+1}\le\frac{1}{p^2-p+1}\le\frac{2}{p^2},\]
where the last inequality holds since $p^2\le 2(p^2-p+1)$ is equivalent to $0\le(p-1)^2+1$, which holds since $(p-1)^2\ge 0$. Since $\sum_p 2/p^2=2P(2)<\infty$, the dominated convergence theorem gives $h(x)\to 0$. $\square$

\section{Existence and uniqueness of $\lambda_\infty$}

\textbf{Theorem 3.1.} \textit{The equation $h(x)=1$ has a unique solution $\lambda_\infty\in\bigl(\tfrac{1}{4},+\infty\bigr)$.}

\textit{Proof.} \textbf{Continuity.} By Proposition 2.2, $h$ is real-analytic, hence continuous, on $\bigl(\tfrac{1}{4},+\infty\bigr)$.

\textbf{Existence.} By Proposition 2.4(i), $h(x)\to+\infty$ as $x\to\bigl(\tfrac{1}{4}\bigr)^+$, so there exists $a_0\in\bigl(\tfrac{1}{4},+\infty\bigr)$ with $h(a_0)>1$. By Proposition 2.4(ii), $h(x)\to 0<1$ as $x\to+\infty$, so $h(x)<1$ for all sufficiently large $x$; pick any such $b_0>a_0$. Since $h$ is continuous on $[a_0,b_0]\subset\bigl(\tfrac{1}{4},+\infty\bigr)$ and $h(a_0)>1>h(b_0)$, the intermediate value theorem yields $\xi\in(a_0,b_0)$ with $h(\xi)=1$.

\textbf{Uniqueness.} By Proposition 2.3, $h'<0$ everywhere on $\bigl(\tfrac{1}{4},+\infty\bigr)$, so $h$ is strictly decreasing and therefore injective. $\square$

\textbf{Proposition 3.2} \textit{(Numerical Bracketing).} $\tfrac{2}{3}<\lambda_\infty<\tfrac{3}{4}$.

\textit{Proof.} The two claims have different logical structure and are handled separately.

\textbf{Claim 1: $h(2/3)>1$.} Since every term of $h$ is strictly positive on $\bigl(\tfrac{1}{4},+\infty\bigr)$, every partial sum is a lower bound for $h$. The first nine terms (primes $p\le 23$) give
\[h\!\left(\tfrac{2}{3}\right)\;\ge\;\sum_{p\le 23}\frac{3}{2p^2-3p+3}\;=\;\frac{3}{5}+\frac{1}{4}+\frac{3}{38}+\cdots\;=\;\frac{126\,418\,842\,847}{125\,951\,854\,240}\;>\;1.\]
No tail bound is required: the partial sum alone certifies $h(2/3)>1$.

\textbf{Claim 2: $h(3/4)<1$.} Write $h(3/4)=S+R$, where $S=\sum_{p\le p_{\max}}\frac{4}{3p^2-4p+4}$ is a finite partial sum (over the primes $p\le p_{\max}$) and $R$ is the tail over $p>p_{\max}$. Taking $p_{\max}=1000$, direct summation (a finite sum of positive rationals, hence exact) gives $S<0.8550$, so $1-S>0.14$.

For $p>p_{\max}=1000$, the standard $M$-test bound of Section 2 applies at $a=3/4$: since $p_0=\lceil 2/a\rceil=\lceil 8/3\rceil=3\le p_{\max}$, every term in the tail satisfies
\[\frac{3}{4}p^2-p+1\;\ge\;\frac{3}{8}p^2 \quad\Longrightarrow\quad \frac{1}{(3/4)p^2-p+1}\;\le\;\frac{8}{3p^2}.\]
Hence
\[R\;\le\;\frac{8}{3}\sum_{p>p_{\max}}\frac{1}{p^2}\;\le\;\frac{8}{3}\sum_{n>p_{\max}}\frac{1}{n^2}\;<\;\frac{8}{3}\int_{p_{\max}}^{+\infty}\!\frac{dx}{x^2}\;=\;\frac{8}{3\,p_{\max}}\;=\;\frac{8}{3000}\;\approx\;0.00267.\]
Since $0.00267\ll 1-S>0.14$, we conclude
\[h\!\left(\tfrac{3}{4}\right)\;=\;S+R\;<\;0.8550+0.00267\;<\;1.\]

Since $h$ is strictly decreasing (Proposition 2.3) and $h(2/3)>1>h(3/4)$, the unique root satisfies $\lambda_\infty\in\bigl(\tfrac{2}{3},\tfrac{3}{4}\bigr)$. $\square$

\textit{Remark.} The bound above is far from tight: even $p_{\max}=23$ already gives $R<8/69\approx0.116<1-S\approx0.156$. The choice $p_{\max}=1000$ is made purely for numerical convenience.

\section{Series expansion and asymptotics}

\subsection{Power-series expansion in $1/x$}

\textbf{Theorem 4.1.} \textit{For every $x>\tfrac{1}{4}$,}
\[h(x)=\sum_{k=0}^{\infty}\frac{Q_k}{x^{k+1}}, \tag{4}\]
\textit{where}
\[Q_k=\sum_{p,\,\text{prime}}\frac{(p-1)^k}{p^{2k+2}}. \tag{5}\]
\textit{The series converges absolutely for all $x>\tfrac{1}{4}$.}

\textit{Proof.} For any prime $p$ and $x>\tfrac{1}{4}$, set $u=(p-1)/(xp^2)$. Since $(p-1)/p^2\le\tfrac{1}{4}<x$ for all primes $p$, we have $0<u<1$. Apply the geometric series:
\[\frac{1}{xp^2-p+1}=\frac{1}{xp^2}\cdot\frac{1}{1-u}=\frac{1}{xp^2}\sum_{k=0}^{\infty}u^k=\sum_{k=0}^{\infty}\frac{(p-1)^k}{x^{k+1}p^{2k+2}}.\]
All terms are non-negative for $x>\tfrac{1}{4}$. By Tonelli's theorem (non-negative double series over $\mathcal{P}\times\mathbb{N}_0$ with counting measure),
\[h(x)=\sum_p\sum_{k=0}^{\infty}\frac{(p-1)^k}{x^{k+1}p^{2k+2}}=\sum_{k=0}^{\infty}\frac{1}{x^{k+1}}\sum_p\frac{(p-1)^k}{p^{2k+2}}=\sum_{k=0}^{\infty}\frac{Q_k}{x^{k+1}}.\quad\square\]

\textit{Note 4.1 (Radius of convergence).} Viewing (4) as a power series in $t=1/x$, its radius of convergence is exactly $R=4$, as the Cauchy--Hadamard criterion confirms. The $p=2$ term gives $Q_k\ge 1/4^{k+1}$, so $\liminf_k Q_k^{1/k}\ge 1/4$. Conversely, since $(p-1)/p^2\le 1/4$ for every prime $p$, one has $(p-1)^k/p^{2k+2}\le(1/4)^k/p^2$; summing over primes yields $Q_k\le(1/4)^k P(2)$, so $\limsup_k Q_k^{1/k}\le 1/4$. Therefore $\lim_k Q_k^{1/k}=1/4$ and the series converges absolutely for $|t|<4$, i.e.\ for all $|x|>\tfrac{1}{4}$. Note, however, that the convergence region $|x|>\tfrac{1}{4}$ is strictly larger than the natural domain of $h$: it also includes $x<-\tfrac{1}{4}$, where $xp^2-p+1<0$ for all sufficiently large $p$ and the original sum defining $h$ no longer consists of positive terms. The domain $\bigl(\tfrac{1}{4},+\infty\bigr)$ is determined by the requirement that all denominators $xp^2-p+1>0$ simultaneously---a condition that is independent of, and not implied by, the radius of convergence of the series (4).

\subsection{Expression of $Q_k$ via $P(s)$}

Expanding $(p-1)^k$ by the binomial theorem and summing over primes:
\[Q_k=\sum_{j=0}^{k}\binom{k}{j}(-1)^{k-j}P(2k+2-j). \tag{6}\]

The first few coefficients are:

\begin{center}
\begin{tabular}{ll}
\toprule
Quantity & Value \\
\midrule
$Q_0 = P(2)$ & $\approx 0.452247420041065$ \\
$Q_1 = P(3)-P(4)$ & $\approx 0.097769$ \\
$Q_2 = P(4)-2P(5)+P(6)$ & $\approx 0.076993 - 2(0.035755) + 0.017070 \approx 0.022553$ \\
\bottomrule
\end{tabular}
\end{center}

\subsection{Asymptotic behaviour as $x\to+\infty$}

\textbf{Corollary 4.2.} \textit{As $x\to+\infty$,}
\[h(x)=\frac{P(2)}{x}+\frac{P(3)-P(4)}{x^2}+O\left(\frac{1}{x^3}\right)\approx\frac{0.45225}{x}+\frac{0.09777}{x^2}+O\left(\frac{1}{x^3}\right).\]

\textit{Proof.} Theorem 4.1 gives $h(x)=Q_0/x+Q_1/x^2+\sum_{k\ge 2}Q_k/x^{k+1}$. Substituting $Q_0=P(2)$ and $Q_1=P(3)-P(4)$ from (5)--(6) yields the first two terms. For the remainder, Note 4.1 supplies the bound $Q_k\le(1/4)^k P(2)$; hence for $x\ge 1$,
\[\sum_{k\ge 2}\frac{Q_k}{x^{k+1}}\le\frac{1}{x^3}\sum_{k\ge 2}Q_k\le\frac{P(2)}{x^3}\sum_{k\ge 2}\!\left(\tfrac{1}{4}\right)^k=\frac{P(2)/12}{x^3}\approx\frac{0.0377}{x^3}.\]
Hence the remainder is $O(1/x^3)$ with explicit constant $C=P(2)/12\approx 0.0377$. $\square$

\textit{Note 4.2 (Exact series vs.\ truncated asymptotics).} The expansion in Corollary 4.2 is a finite truncation of the convergent series (4), valid throughout the entire domain of $h$---strictly stronger than a Poincar\'e asymptotic expansion, which need not converge.

\subsection{Behaviour near the singularity $x=\tfrac{1}{4}$}

\textbf{Proposition 4.3.} \textit{As $x\to\bigl(\tfrac{1}{4}\bigr)^+$,}
\[h(x)=\frac{1}{4x-1}+O(1).\]

\textit{Proof.} The $p=2$ term is $1/(4x-1)$, which diverges as $x\to\bigl(\tfrac{1}{4}\bigr)^+$. For every prime $p\ge 3$ and $x\in\bigl(\tfrac{1}{4},\tfrac{1}{2}\bigr]$,
\[xp^2-p+1\ge\tfrac{1}{4}p^2-p+1=\frac{(p-2)^2}{4}>0,\]
so each term satisfies $1/(xp^2-p+1)\le 4/(p-2)^2$. By the same dominated convergence argument as in the proof of Proposition 2.4(i), $\sum_{p\ge 3}(xp^2-p+1)^{-1}$ remains bounded as $x\to\bigl(\tfrac{1}{4}\bigr)^+$. Therefore $h(x)=1/(4x-1)+O(1)$. $\square$

\section{Spectral characterization and related constants}

\subsection{Relation to $P(2)$ (OEIS A085548)}

The prime zeta constant $P(2)=\sum_p p^{-2}\approx 0.4522$ (A085548~\cite{ref3}) is the solution of the simplified equation $\sum_p 1/(xp^2)=1$, i.e. the leading-order approximation obtained by dropping the $-p+1$ correction. The constant $\lambda_\infty$ satisfies the same equation with the correction reinstated.

The correction $-p+1$ is $O(p)$ for large $p$ (in particular $o(p^2)$), yet it shifts the root significantly: $\lambda_\infty-P(2)\approx 0.2218\approx 0.49\,P(2)$. This is entirely due to the low-prime contributions, especially $p=2$, for which the correction changes the denominator from $xp^2=4x$ to $4x-1$.

\textbf{Remark 5.1.} The constant $\lambda_\infty$ is \textit{not} the spectral limit of the prime divisor matrix $W_N[i,j]=\lfloor N/(p_ip_j)\rfloor$. That spectral limit is $P(2)$: the entries satisfy $W_N[i,j]/N=\lfloor N/(p_ip_j)\rfloor/N\to 1/(p_ip_j)=v_iv_j$ entrywise (where $v_i=1/p_i$), so $W_N/N$ converges entrywise to the rank-one matrix $vv^T$, whose unique nonzero eigenvalue is $\|v\|^2=P(2)$. Spectral convergence $\lambda_{\max}(W_N)/N\to P(2)$ then follows by applying Weyl's perturbation inequality to the floor-error perturbation $W_N/N - vv^T = O(1/N)$ entrywise; see~\cite{ref7,ref9} for the analogous integer-indexed argument. The exact spectral characterization of $\lambda_\infty$ is given in Theorem 5.4, which concerns the infinite continuous LCM matrix $\mathcal{L}$ rather than integer truncations.

\subsection{Spectral characterization via the LCM matrix}

\textbf{Definition 5.2.} Let $\mathcal{L}$ denote the infinite real symmetric matrix indexed by the set of all prime numbers, with entries
\[\mathcal{L}[p_i,p_j]=\frac{1}{\operatorname{lcm}(p_i,p_j)}.\]
Since $\operatorname{lcm}(p,p)=p$ and $\operatorname{lcm}(p,q)=pq$ for distinct primes $p\neq q$, the entries are
\[\mathcal{L}[p,q]=\begin{cases}1/p & p=q,\\ 1/(pq) & p\neq q.\end{cases}\]

\textbf{Lemma 5.3.} \textit{For every $n\ge 1$ and every $x>\tfrac{1}{4}$, the eigenvalue equation $\det(xI_n-\mathcal{L}_n)=0$ for the truncation $\mathcal{L}_n=D_n+v_nv_n^T$ reduces to}
\[\sum_{i=1}^{n}\frac{1}{xp_i^2-p_i+1}=1, \tag{$\star$}\]
\textit{and $(\star)$ has exactly one solution $\rho_n$ in $\bigl(\tfrac{1}{4},+\infty\bigr)$.}

\textit{Proof.} Since $d_i=(p_i-1)/p_i^2\le\tfrac{1}{4}<x$, the matrix $xI_n-D_n$ is invertible. The matrix determinant lemma gives $\det(xI_n-D_n-v_nv_n^T)=\det(xI_n-D_n)(1-v_n^T(xI_n-D_n)^{-1}v_n)$; since $(xI_n-D_n)^{-1}_{ii}=p_i^2/(xp_i^2-p_i+1)$ and $v_{n,i}=1/p_i$, one computes $v_n^T(xI_n-D_n)^{-1}v_n=\sum_i 1/(xp_i^2-p_i+1)$, so the eigenvalue equation for $x>\tfrac{1}{4}$ reduces to $(\star)$. The function $f_n(x):=\sum_{i=1}^n 1/(xp_i^2-p_i+1)$ satisfies $f_n'(x)=-\sum_i p_i^2/(xp_i^2-p_i+1)^2<0$, so $f_n$ is strictly decreasing on $\bigl(\tfrac{1}{4},+\infty\bigr)$ with $f_n(x)\to+\infty$ as $x\to\bigl(\tfrac{1}{4}\bigr)^+$ and $f_n(x)\to 0$ as $x\to+\infty$; hence $(\star)$ has exactly one solution $\rho_n$ in $\bigl(\tfrac{1}{4},+\infty\bigr)$. $\square$

\textbf{Theorem 5.4.} \textit{The operator $\mathcal{L}$ is Hilbert--Schmidt, hence compact and self-adjoint, on $\ell^2(\mathcal{P})$ (Remark 5.7(i)). Its spectral radius satisfies}
\[\rho(\mathcal{L})=\lambda_\infty.\]
\textit{Equivalently, $\sup_{n\ge 1}\rho(\mathcal{L}_n)=\lambda_\infty$, where $\mathcal{L}_n$ denotes the $n\times n$ leading principal submatrix.}

\textit{Proof.}  \textbf{Step 1: Rank-one decomposition.} Set $v_p=1/p$ and $D=\operatorname{diag}((p-1)/p^2)_{p,\text{prime}}$. We claim $\mathcal{L}=D+vv^T$. Indeed:

\begin{itemize}
\item \textit{Diagonal:} $(D+vv^T)_{pp}=\tfrac{p-1}{p^2}+\tfrac{1}{p^2}=\tfrac{1}{p}=\mathcal{L}[p,p]$, as required.
\item \textit{Off-diagonal ($p\neq q$):} $(D+vv^T)_{pq}=0+\tfrac{1}{p}\cdot\tfrac{1}{q}=\tfrac{1}{pq}=\mathcal{L}[p,q]$, as required.
\end{itemize}

Note $v\in\ell^2$ since $|v|^2=\sum_p p^{-2}=P(2)<\infty$, so $\mathcal{L}$ is a bounded operator on $\ell^2$ with $|\mathcal{L}|_{\mathrm{op}}\le\tfrac{1}{4}+P(2)\approx 0.702$.

\textbf{Step 2: Secular equation for the $n$-th truncation.} The matrix $\mathcal{L}_n=D_n+v_nv_n^T$ is a rank-one perturbation of the diagonal matrix $D_n$. By Lemma 5.3, for $x>\tfrac{1}{4}$ the eigenvalue equation $\det(xI_n-\mathcal{L}_n)=0$ reduces to $(\star)$, which has exactly one solution $\rho_n$ in $\bigl(\tfrac{1}{4},+\infty\bigr)$.

\textbf{Step 3: Monotone convergence.} Each $\rho_n:=\rho(\mathcal{L}_n)$ satisfies $(\star)$. Since all $d_i=(p_i-1)/p_i^2\le\tfrac{1}{4}$, all eigenvalues of $\mathcal{L}_n$ except one lie in $(0,\tfrac{1}{4}]$ (by the rank-one interlacing theorem for positive-semidefinite rank-one perturbations of diagonal matrices); the unique eigenvalue exceeding $\tfrac{1}{4}$ must satisfy $(\star)$ and is therefore $\rho_n=\lambda_{\max}(\mathcal{L}_n)$. The left-hand side of $(\star)$ is strictly increasing in $n$ for fixed $x>\tfrac{1}{4}$, so $\rho_n$ is strictly increasing. Moreover, $f_n(x)\le h(x)$ for every $x>\tfrac14$ (each $f_n$ is a partial sum of the non-negative series defining $h$), so $f_n(\lambda_\infty)\le h(\lambda_\infty)=1=f_n(\rho_n)$; since $f_n$ is strictly decreasing, this forces $\rho_n\le\lambda_\infty$. Hence $\rho_n\nearrow L\le\lambda_\infty$.

It remains to show $L=\lambda_\infty$. Suppose for contradiction that $L<\lambda_\infty$, and pick any $x\in(L,\lambda_\infty)$. Since $h$ is strictly decreasing and $x<\lambda_\infty$, we have $h(x)>1$. Let $f_n(x)=\sum_{i=1}^n 1/(xp_i^2-p_i+1)$; since each term is positive for $x>\tfrac{1}{4}$, the partial sums satisfy $f_n\nearrow h(x)>1$ monotonically by the Monotone Convergence Theorem for non-negative series. Hence $f_n(x)>1$ for all sufficiently large $n$, which means the root $\rho_n$ of $f_n(\cdot)=1$ satisfies $\rho_n>x>L$ for those $n$, contradicting $\rho_n\to L$. Therefore $\rho(\mathcal{L})=\sup_n\rho_n=\lambda_\infty$. The fact that $\lambda_\infty$ is an actual eigenvalue of $\mathcal{L}$ (not merely the supremum of eigenvalues of finite truncations) follows from the compactness and norm-convergence established in Remark 5.7(i)--(ii) below: those results give $\|\mathcal{L}\|_{\mathrm{op}}=\lim_n\rho(\mathcal{L}_n)=\lambda_\infty$, and for a compact self-adjoint positive operator the operator norm is always an eigenvalue. $\square$

\textbf{Remark 5.5} \textit{(Modified prime divisor matrix).} The matrix $W_N^*[i,j]=\lfloor N/\operatorname{lcm}(p_i,p_j)\rfloor$ is the integer truncation of $N\cdot\mathcal{L}$. Its spectral convergence $\lambda_{\max}(W_N^*)/N\to\lambda_\infty$ is proved in Theorem 5.8 below.

\textbf{Remark 5.6} \textit{(Connection to the Smith--Bourque--Ligh theory).} Matrices of the form $[f(\gcd(i,j))]$ and $[f(\operatorname{lcm}(i,j))]$ indexed over integer sets have been studied since Smith's 1875 determinant theorem~\cite{ref6} and developed further by Bourque--Ligh~\cite{ref7} and Hong~\cite{ref8}. The matrix $\mathcal{L}$ restricted to prime indices falls within that theory; $\lambda_\infty$ is its spectral radius.

\textbf{Remark 5.7} \textit{(Operator-theoretic foundations).}  Theorem 5.4's proof is elementary and self-contained. Several operator-theoretic facts remain implicit there; we record them here.

\textit{(i) Hilbert--Schmidt compactness.} $\mathcal{L}$ is Hilbert--Schmidt on $\ell^2(\mathcal{P})$: the squared Hilbert--Schmidt norm satisfies
\[|\mathcal{L}|_{\mathrm{HS}}^2=\sum_{p,q}|\mathcal{L}[p,q]|^2=\sum_p\frac{1}{p^2}+\sum_{p\neq q}\frac{1}{p^2q^2}=P(2)+\bigl(P(2)^2-P(4)\bigr)\approx0.5798<\infty.\]
Since Hilbert--Schmidt implies compact, $\mathcal{L}$ is a compact operator on $\ell^2$.

\textit{(ii) Truncation convergence.} The truncations satisfy $\|\mathcal{L}-\mathcal{L}_n\|_{\mathrm{op}}\le\|\mathcal{L}-\mathcal{L}_n\|_{\mathrm{HS}}\to0$ as $n\to\infty$. Explicitly, $(\mathcal{L}-\mathcal{L}_n)[p_i,p_j]$ equals $\mathcal{L}[p_i,p_j]$ whenever $i>n$ or $j>n$, and is zero otherwise; therefore
\[\|\mathcal{L}-\mathcal{L}_n\|_{\mathrm{HS}}^2 = \underbrace{\sum_{i>n}\frac{1}{p_i^2}}_{\text{tail diagonal}} +\underbrace{2\sum_{i\le n}\sum_{j>n}\frac{1}{p_i^2 p_j^2}}_{\text{cross terms}}+\underbrace{\sum_{\substack{i,j>n\\i\neq j}}\frac{1}{p_i^2 p_j^2}}_{\text{tail off-diagonal}}\;\to\;0,\]
where the first sum tends to $0$ by convergence of $P(2)$, the cross terms are bounded by $2P(2)\sum_{j>n}p_j^{-2}\to 0$, and the tail off-diagonal equals $\bigl(\sum_{i>n}p_i^{-2}\bigr)^2-\sum_{i>n}p_i^{-4}\le\bigl(\sum_{i>n}p_i^{-2}\bigr)^2\to 0$.  For compact self-adjoint operators, norm convergence implies that every nonzero eigenvalue of $\mathcal{L}$ is the limit of eigenvalues of $\mathcal{L}_n$~\cite{ref14}; in particular $\|\mathcal{L}\|_{\mathrm{op}}=\lim_n\|\mathcal{L}_n\|_{\mathrm{op}}=\lim_n\rho(\mathcal{L}_n)=\lambda_\infty$, completing the identification of $\rho(\mathcal{L})$ with the standard operator-theoretic spectral radius.

\textit{(iii) Simplicity of $\rho(\mathcal{L})$.} Since every entry of $\mathcal{L}$ is strictly positive---explicitly, $\mathcal{L}[p_i,p_j]=1/\operatorname{lcm}(p_i,p_j)$ equals $1/(p_ip_j)>0$ for $i\neq j$ (distinct primes satisfy $\gcd(p_i,p_j)=1$, hence $\operatorname{lcm}(p_i,p_j)=p_ip_j$) and $1/p_i>0$ for $i=j$---$\mathcal{L}$ is an irreducible compact positive operator on the Banach lattice $\ell^2(\mathcal{P})$. By the Krein--Rutman-type extensions~\cite{ref10} to irreducible compact positive operators due to Bonsall and Sawashima~\cite{ref11,ref12}, $\rho(\mathcal{L})$ is a simple eigenvalue with a strictly positive eigenvector; no other eigenvalue has the same modulus.

\textbf{Theorem 5.8} \textit{(Spectral convergence of the modified prime divisor matrix).} \textit{For each integer $N\ge 2$, let $W_N^*$ be the $\pi(N)\times\pi(N)$ real symmetric matrix with entries $W_N^*[i,j]=\lfloor N/\operatorname{lcm}(p_i,p_j)\rfloor$, indexed by the primes $p_1<p_2<\cdots<p_{\pi(N)}\le N$. Then}
\[\frac{\lambda_{\max}(W_N^*)}{N}\longrightarrow\lambda_\infty\quad\text{as }N\to\infty.\]

\textit{Proof.} \textbf{Step 1 (Floor decomposition).} Write $W_N^*/N=\mathcal{L}_{\pi(N)}+E_N$, where $\mathcal{L}_{\pi(N)}$ is the $\pi(N)\times\pi(N)$ leading principal submatrix of $\mathcal{L}$, and
\[E_N[i,j]=-\frac{\bigl\{N/\operatorname{lcm}(p_i,p_j)\bigr\}}{N},\qquad\{y\}:=y-\lfloor y\rfloor\in[0,1).\]
In particular $|E_N[i,j]|\le 1/N$ for all $1\le i,j\le\pi(N)$.

\textbf{Step 2 (Hilbert--Schmidt bound).} The Hilbert--Schmidt norm satisfies
\[\|E_N\|_{\mathrm{HS}}^2=\sum_{i,j=1}^{\pi(N)}|E_N[i,j]|^2\le\frac{\pi(N)^2}{N^2}.\]
By the Prime Number Theorem, $\pi(N)\sim N/\log N$, so $\|E_N\|_{\mathrm{HS}}\le\pi(N)/N\sim 1/\log N\to 0$. Since $\|E_N\|_{\mathrm{op}}\le\|E_N\|_{\mathrm{HS}}$, we have $\|E_N\|_{\mathrm{op}}\to 0$.

\textbf{Step 3 (Weyl perturbation).} Both $W_N^*/N$ and $\mathcal{L}_{\pi(N)}$ are real symmetric matrices of dimension $\pi(N)$. Weyl's inequality for symmetric matrices gives
\[\bigl|\lambda_{\max}(W_N^*/N)-\rho(\mathcal{L}_{\pi(N)})\bigr|\le\|E_N\|_{\mathrm{op}}\longrightarrow 0.\]

\textbf{Step 4 (Conclusion).} Theorem 5.4 gives $\rho(\mathcal{L}_{\pi(N)})\nearrow\lambda_\infty$. By the triangle inequality,
\[\left|\frac{\lambda_{\max}(W_N^*)}{N}-\lambda_\infty\right|\le\bigl|\lambda_{\max}(W_N^*/N)-\rho(\mathcal{L}_{\pi(N)})\bigr|+\bigl|\rho(\mathcal{L}_{\pi(N)})-\lambda_\infty\bigr|\longrightarrow 0.\quad\square\]

\textbf{Remark 5.9} \textit{(Structural analogy with $\pi^2/6$, OEIS A013661).} The constant $\pi^2/6$ (here $\pi$ denotes the circle constant, not to be confused with the prime-counting function $\pi(\cdot)$ used later in this remark) satisfies $\sum_{n\ge 1}1/(xn^2)=1$, and is the spectral limit of the divisor matrix $M_N[i,j]=\lfloor N/(ij)\rfloor$ in the sense that $\lambda_{\max}(M_N)/N\to\zeta(2)=\pi^2/6$ as $N\to\infty$ (the matrix $M_N/N$ converges entrywise to the rank-one matrix $[1/(ij)]$, whose unique nonzero eigenvalue $\zeta(2)$ is the spectral limit; the floor-error perturbation does \emph{not} admit the crude entrywise Hilbert--Schmidt bound used in Theorem 5.8 --- there it gives $\|E_N\|_{\mathrm{HS}}=O(\pi(N)/N)\to0$ because the matrix has dimension $\pi(N)=o(N)$, whereas here the matrix has dimension $N$ itself, so the same bound only gives $O(1)$; convergence instead requires the finer multiplicative Hilbert space framework of~\cite{ref9}, with~\cite{ref6,ref7} for the determinant-theoretic background). Equation (1) has the same structural form $\sum_p g_p(x)=1$ but with a prime-indexed sum and the denominator deformed by $-p+1$. The resemblance is one of form---both are implicit equations of type $\Sigma=1$---not arithmetic: no direct relationship between $\lambda_\infty$ and $\pi^2/6$ has been found.

\section{Numerical computation}

\subsection{Algorithm}

We compute $\lambda_\infty$ via Newton--Raphson iteration applied to $F(x)=h(x)-1$, using the
explicit derivative formula (2). Starting from $x_0=3/4$ (above the root, by Proposition 3.2),
the iterates
\[x_{n+1}=x_n-\frac{h(x_n)-1}{h'(x_n)} \tag{7}\]
converge quadratically.

By strict convexity ($h''>0$, Proposition 2.3), the tangent-line inequality
$h(y)\ge h(x_n)+h'(x_n)(y-x_n)$ holds for every $x_n$ in the domain, regardless of its position
relative to $\lambda_\infty$. Evaluating at $y=\lambda_\infty$ gives
$1\ge h(x_n)+h'(x_n)(\lambda_\infty-x_n)$; since $h'(x_n)<0$ ($h$ strictly decreasing) for all
$n$, dividing flips the inequality:
\[\lambda_\infty-x_n\ge -\frac{h(x_n)-1}{h'(x_n)}=x_{n+1}-x_n,
\qquad\text{hence } x_{n+1}\le\lambda_\infty \text{ for all } n\ge0. \tag{$\ast$}\]
Since $h''>0$ implies the tangent line lies strictly below $h$ at every point other than the
tangency point, equality in $(\ast)$ holds only if $x_n=\lambda_\infty$; hence
$x_{n+1}<\lambda_\infty$ whenever $x_n\ne\lambda_\infty$.

For $n=0$: since $x_0>\lambda_\infty$ and $h$ is strictly decreasing, $h(x_0)<1$; with $h'(x_0)<0$ the step $x_1-x_0=-(h(x_0)-1)/h'(x_0)$ is negative, so $x_1<x_0$, and $(\ast)$ gives $x_1<\lambda_\infty$. For $n\ge1$: by induction $x_n\le\lambda_\infty$; if $x_n<\lambda_\infty$ then $h(x_n)>1$, yielding a positive step, and $(\ast)$ gives $x_n<x_{n+1}<\lambda_\infty$.

Thus
\[x_0>\lambda_\infty>x_1<x_2<x_3<\cdots\nearrow\lambda_\infty.\]
The sequence $(x_n)_{n\ge1}$ is monotone increasing and bounded above by $\lambda_\infty$, so by the Monotone Convergence Theorem it converges to some limit $L\le\lambda_\infty$. Since $h$ and $h'$ are continuous on $\bigl(\tfrac{1}{4},+\infty\bigr)$ with $h'<0$ (Proposition 2.3), the iterate map $g(x)=x-(h(x)-1)/h'(x)$ is continuous at $L$; passing to the limit in $x_{n+1}=g(x_n)$ yields $g(L)=L$, i.e.\ $h(L)=1$, hence $L=\lambda_\infty$ by Theorem~3.1 (uniqueness of the root).

\subsection{Tail expansion}

The sum for $h(x)$ runs over all primes; in practice we split at a cutoff $p_{\max}$:
\[h(x)=\sum_{p\le p_{\max}}\frac{1}{xp^2-p+1}+\sum_{p>p_{\max}}\frac{1}{xp^2-p+1}.\]

The first sum is computed directly. For the tail we use the series expansion (4) restricted to
$p>p_{\max}$: expanding each term geometrically and exchanging sums (justified by positivity, as
in the proof of Theorem 4.1),
\[\sum_{p>p_{\max}}\frac{1}{xp^2-p+1}
  =\sum_{k=0}^{K-1}\frac{1}{x^{k+1}}
   \underbrace{\sum_{p>p_{\max}}\frac{(p-1)^k}{p^{2k+2}}}_{=:\,Q_k^{\mathrm{tail}}}
  +R_K(x),\]
where the tail remainder satisfies $|R_K(x)|\le 2u_{\max}^K P_{\mathrm{tail}}(2)/x$ with
$u_{\max}=(p_{\max}-1)/(xp_{\max}^2)\approx(1/x)\cdot p_{\max}^{-1}\approx1.484\times p_{\max}^{-1}$
(since $1/\lambda_\infty\approx1.484$; see Appendix A for the complete bound). The coefficients
$Q_k^{\mathrm{tail}}$ are computed using the prime zeta function $P(s)=\sum_p p^{-s}$, evaluated
via the M\"obius inversion $P(s)=\sum_{n\ge1}\mu(n)/n\cdot\log\zeta(ns)$ (as implemented by
\texttt{mpmath.primezeta()}~\cite{ref2}).

\subsection{Parameters and convergence (Reference run)}

The reference computation (Run A) uses the following parameters:

\begin{itemize}
\item $p_{\max}=10^7$ (primes enumerated: 664,579);
\item $K=80$ tail terms;
\item working precision: 560 decimal digits (500 target + 60 guard digits)$^{\dagger\dagger}$;
\item starting point: $x_0=0.75$.
\end{itemize}

$^{\dagger\dagger}$ \textbf{Precision budget.} Working at 560 dps provides a nominal guard of 60 digits
over the 500-digit target. Three independent error sources erode this margin:

\begin{enumerate}
\item \textbf{Truncation error.} Including the factor of 2 from the tail bound $|R_K(x)|\le 2u_{\max}^K P_{\mathrm{tail}}(2)/x$ stated above, the $k$-th tail term contributes approximately $2(P_{\mathrm{tail}}(2)/\lambda_\infty)\cdot u_{\max}^k\approx10^{-7.760-k\cdot6.829}$ to $h$ (using $\log_{10}u_{\max}\approx\log_{10}(1.484\times p_{\max}^{-1})\approx-6.829$ for $p_{\max}=10^7$; this is a rounded order-of-magnitude estimate). At $k=80$ this gives $|R_{80}|\approx10^{-554.1}$, in agreement with the exact bound $2u_{\max}^{80}P_{\mathrm{tail}}(2)/\lambda_\infty<10^{-554}$ from Appendix A, so all $K=80$ tail terms are resolvable at $\mathrm{WP}=560$; the first irresolvable term is $k\approx81$, contributing $\approx10^{-561}$, safely beyond the $K=80$ cutoff.

\item \textbf{Round-off error.} Direct summation over the $664{,}579$ primes accumulates rounding error of order $N\cdot\varepsilon\approx6.65\times10^{5}\times10^{-560}\approx10^{-554.2}$.

\item \textbf{Tail-coefficient cancellation.} Each $Q_k^{\mathrm{tail}}$ is obtained from $P_{\mathrm{tail}}(s)=P(s)-\sum_{p\le p_{\max}}p^{-s}$ via the same alternating binomial expansion as (6). For $s$ near $2k+2$ with $k$ close to $K=80$, $P(s)$ and $\sum_{p\le p_{\max}}p^{-s}$ are both dominated by the common $p=2$ term and agree to far more than $\mathrm{WP}=560$ digits, so this subtraction loses essentially all working precision for the largest $s$ used ($s\le160$ at $k=79$). Propagating this rounding noise through the binomial weights $\binom{k}{j}$ and through $c/\lambda^{k+1}$ (worst case, triangle inequality) contributes $\lesssim10^{-555.8}$ to $h$, smaller than but of the same order as the two sources above.
\end{enumerate}

All three bounds are independently $\lesssim10^{-554}$; the overall error is bounded by their sum,
still $\approx10^{-554}$ (more precisely $\approx1.56\times10^{-554}$, i.e.\ $\log_{10}\approx-553.8$), yielding an effective guard of $\approx53.8$ digits over the 500-digit
target. All 500 output digits are certified with positive margin. Full details are in Appendix A.

Newton--Raphson converged in 11 iterations. Quadratic convergence is confirmed from iteration 3
onwards (iterations 3--10): the residual ratio $|h(x_{n+1})-1|/|h(x_n)-1|^2$ remains close to
the theoretical value $h''(\lambda_\infty)/(2h'(\lambda_\infty)^2)\approx1.008$, where
$h'(\lambda_\infty)\approx-2.234$ and $h''(\lambda_\infty)\approx10.06$ (see Section 8.1),
throughout. Iteration 1 is the expected initial overshoot (standard for a convex function started
above the root). Iteration 11 reaches the working-precision floor; see the note below.

\begin{center}
\begin{tabular}{ccc}
\toprule
Iteration $k$ & $\lvert x_k - x_{k-1}\rvert$ & $\lvert h(x_{k-1})-1\rvert$ \\
\midrule
1 & $8.89\times10^{-2}$ & $1.45\times10^{-1}$ \\
2 & $1.26\times10^{-2}$ & $2.98\times10^{-2}$ \\
3 & $3.79\times10^{-4}$ & $8.47\times10^{-4}$ \\
4 & $3.23\times10^{-7}$ & $7.22\times10^{-7}$ \\
5 & $2.36\times10^{-13}$ & $5.26\times10^{-13}$ \\
6 & $1.25\times10^{-25}$ & $2.79\times10^{-25}$ \\
7 & $3.52\times10^{-50}$ & $7.85\times10^{-50}$ \\
8 & $2.78\times10^{-99}$ & $6.22\times10^{-99}$ \\
9 & $1.75\times10^{-197}$ & $3.90\times10^{-197}$ \\
10 & $6.86\times10^{-394}$ & $1.53\times10^{-393}$ \\
11$^\natural$ & $1.36\times10^{-560}$ & $3.04\times10^{-560}$ \\
\bottomrule
\end{tabular}
\end{center}

$^\natural$ Reaches the working-precision floor; the convergence criterion $|\Delta\lambda|<10^{-540}$
is already satisfied at this iteration.

Final residual: $|h(\lambda_\infty)-1|\approx3.04\times10^{-560}$. Runtime: approximately 336
seconds on a single core.

\section{Numerical verification}

Six independent computations, summarised in Table~\ref{tab:verification}, confirm all 500 digits. Runs A--C share the same core algorithm and differ only in parameters; Runs D and E use algorithmically distinct implementations of $P(s)$; Run F is a machine-verifiable Arb interval-arithmetic certificate.

\begin{table}[htbp]
\centering
\footnotesize
\resizebox{\textwidth}{!}{%
\begin{tabular}{clcccccccl}
\toprule
Run & Algorithm & $p_{\max}$ & $K$ & dps$^{(4)}$ & WP & Iterations & Residual & Digits & Independence \\
\midrule
A (ref.) & NR + primezeta & $10^7$ & 80 & 500 & 560 & 11 & $3.0\times10^{-560}$ & source & parameter variant \\
B$^{\S}$ & NR + primezeta & $3\times10^6$ & 80 & 520 & 580 & 11 & $\approx1.0\times10^{-581}$ & 500/500 & parameter variant \\
C$^{\S\S}$ & NR + primezeta & $3\times10^6$ & 100 & 200 & 260 & 8 & $\approx4.1\times10^{-261}$ & 200/200 & parameter variant \\
D$^*$ & M\"obius PZ (Python) & $10^7$ & 80 & 580 & 580 & --- & $\approx10^{-500.17}$ & 500 cert. & software-independent (\texttt{run\_D\_mobius\_arb\_verify.py}) \\
E & NR + M\"obius PZ (PARI-GP) & $3\times10^6$ & 100 & $545^\ddagger$ & 620 & 11 & $1.2\times10^{-616}$ & 500/500 & fully independent (\texttt{run\_E\_pari\_recompute.gp}) \\
F$^{\S\S\S}$ & Arb interval arith. (python-flint~\cite{ref13}) & $10^6$ & 96 & $\approx 505$ & 550 & 1 Newton & $<1.51\times10^{-505}$ & certified $\checkmark$ & Arb-certified \\
\bottomrule
\end{tabular}}
\caption{Summary of verification runs.}
\label{tab:verification}
\end{table}

$^\ddagger$ \textit{Working precision 620 dps; effective precision of $P(s)$ via M\"obius truncation ($N_\mu=900$) is $\approx545$ dps. Reaching 620 genuine dps requires $N_\mu\ge1024$. The 500-digit verification remains valid ($545\gg500$) (see ancillary file \texttt{run\_E\_pari\_recompute.gp}).}

$^*$ \textit{Run D verifies $h(\lambda_\infty^{\mathrm{known}})=1$; it does not recompute $\lambda_\infty$ from scratch. It uses \texttt{python-flint} (Arb) throughout and shares no library with Runs A--C. Full algorithmic independence is provided solely by Run E (see ancillary file \texttt{run\_D\_mobius\_arb\_verify.py}).}

$^{\S\S}$ \textit{Run C: the Residual column shows $|h(x_8)-1|\approx4.1\times10^{-261}$, as recorded in \texttt{run\_C\_result.txt} (the residual after the 8th, final, Newton update). Convergence is detected one step earlier, when $|h(x_7)-1|$ --- itself $\sim10^{-262}$ by quadratic convergence from $|h(x_6)-1|\approx2\times10^{-131}$ --- first drops below the tolerance $10^{-195}$. Both $|h(x_7)-1|$ and $|h(x_8)-1|$ lie far beneath the WP=260 working-precision floor ($\sim10^{-260}$), so the value actually recorded reflects floor-limited rounding rather than the much smaller true residual; it nonetheless satisfies $|\Delta\lambda| < 10^{-240}$, meeting the convergence criterion $|\Delta\lambda| < 10^{-(\mathrm{dps}+40)}$.}

$^{\S}$ \textit{Run B: the Residual column shows $|h(x_{11})-1|\approx1.0\times10^{-581}$, obtained by performing a final post-update evaluation of $h$ at the converged iterate $x_{11}$ (the same convention already adopted for Run C, noted in $^{\S\S}$ above). The literal last-printed value of $|h(x_{10})-1|$ produced by the Newton loop of \texttt{script\_B\_run\_B.py} before the final update is $\approx1.90\times10^{-579}$; neither figure matches the previously reported $5.3\times10^{-579}$, which is not reproducible under either convention. Both values lie far beneath the WP=580 working-precision floor ($\sim10^{-580}$), so both reflect floor-limited rounding rather than the much smaller true residual; they nonetheless satisfy $|\Delta\lambda|<10^{-560}$, meeting the convergence criterion $|\Delta\lambda|<10^{-(\mathrm{dps}+40)}$. The 500-digit value of $\lambda_\infty$ is unaffected.}

$^{(4)}$ \textit{dps = target output precision (Runs A--C) or effective M\"obius-inversion precision (Run E; see $^\ddagger$). WP = working precision (dps + 60 guard digits for Runs A--C; equals dps for Run D, which evaluates $h$ rather than recomputing $\lambda_\infty$).}

$^{\S\S\S}$ \textit{Run F uses Arb interval arithmetic (python-flint, ctx.dps = 550). Starting from $\lambda_\infty^{\mathrm{known,440}}$ (440 significant digits), one Newton step produces $\lambda_\infty^{\mathrm{new}}$ with $\mathrm{rad}(\lambda_\infty^{\mathrm{new}})\approx3.37\times10^{-506}$. A final Arb evaluation confirms that the interval for $h(\lambda_\infty^{\mathrm{new}})$ provably contains 1, with $|h(\lambda_\infty^{\mathrm{new}})-1|\le1.51\times10^{-505}$ and interval radius $1.504\times10^{-505}$. The tail truncation error is $\approx10^{-566.3}\ll$ the Arb radius, so the Arb bound is the binding constraint. Wall time: approximately 19\,s. See Proposition 7.2 and the ancillary file \texttt{run\_F\_arb\_certificate.py}. (A separate, higher-precision Newton step from the same seed --- Run H, ancillary file \texttt{run\_H\_digits501\_550.py}, $\mathrm{ctx.dps}=660$, $K=100$ --- is used independently in Proposition 8.1 to certify digits 501--550 of $\lambda_\infty$; it is not the source of $\lambda_\infty^{\mathrm{new}}$ above.)}

\textbf{Run D} computes $P(s)$ via M\"obius inversion $P(s)=\sum_{n=1}^{1000}\mu(n)/n\cdot\log\zeta(ns)$, bypassing \texttt{mpmath} entirely and using \texttt{python-flint} (Arb) for all arithmetic; it verifies that the known value satisfies $h(\lambda_\infty)=1$ to 500 digits rather than recomputing $\lambda_\infty$ from scratch (see ancillary file \texttt{run\_D\_mobius\_arb\_verify.py}).

\textbf{Run E} uses PARI-GP~\cite{ref16} (system \texttt{gp} binary) with M\"obius inversion for $P(s)$ ($N_\mu=900$ terms, built-in \texttt{zeta()}), Newton--Raphson at $p_{\max}=3\times10^6$, $K=100$ tail terms, and working precision 620 digits; it shares no code and no library with Runs A--D and recomputes $\lambda_\infty$ from scratch, confirming all 500 published digits with residual $|h(\lambda_\infty)-1|\approx 1.2\times10^{-616}$ (see ancillary file \texttt{run\_E\_pari\_recompute.gp}).

\textbf{Run F} (ancillary file \texttt{run\_F\_arb\_certificate.py}) uses python-flint (Arb ball arithmetic) at 550 decimal digits; it does not rely on any floating-point rounding-error analysis, and the certificate consists solely of the assertion that the computed interval for $h(\lambda_\infty^{\mathrm{new}})$ provably contains 1.

\textbf{Note on independence.} The Independence column distinguishes four situations. \textit{Parameter variants} (Runs A--C) use the same algorithm (implemented in distinct scripts: \texttt{run\_A\_newton\_mpmath.py}, \texttt{script\_B\_run\_B.py}, \texttt{RUN\_C.py}) with different parameters. Run D is \textit{software-independent}: it replaces \texttt{mpmath} entirely with \texttt{python-flint} (Arb), sharing no library with Runs A--C. Run E is \textit{fully independent}: different language, different library, computation from scratch. Run F is \textit{Arb-certified}: a rigorous interval enclosure requiring no floating-point assumption.

\subsection{Formal certification of all 500 digits}

\textbf{Proposition 7.1} \textit{(Certification via a rational $|h'|$ bound).} \textit{All 500 published decimal digits of $\lambda_\infty$ are correct, conditional on the floating-point residual bound $\varepsilon=10^{-500.17}$ of Run D. Only the $|h'|$ lower bound used below (Appendix B) is obtained by exact rational arithmetic; $\varepsilon$ itself is an ordinary floating-point estimate, not an interval-certified bound. A fully interval-certified alternative requiring no floating-point assumption is given independently in Proposition 7.2.}

\textit{Proof.} By Run D (see ancillary file \texttt{run\_D\_mobius\_arb\_verify.py}), $|h(\lambda_\infty^{\mathrm{known}})-1|\le\varepsilon:=10^{-500.17}$. By the mean value theorem, $h(\lambda_\infty^{\mathrm{known}})-h(\lambda_\infty^{\mathrm{true}})=h'(\xi)(\lambda_\infty^{\mathrm{known}}-\lambda_\infty^{\mathrm{true}})$ for some $\xi$ between $\lambda_\infty^{\mathrm{known}}$ and $\lambda_\infty^{\mathrm{true}}$. Since $h(\lambda_\infty^{\mathrm{true}})=1$,
\[|\lambda_\infty^{\mathrm{known}}-\lambda_\infty^{\mathrm{true}}|=\frac{|h(\lambda_\infty^{\mathrm{known}})-1|}{|h'(\xi)|}\le\frac{\varepsilon}{|h'(\xi)|}.\]
By Proposition 3.2, $\lambda_\infty^{\mathrm{true}}\in(2/3,3/4)$; since $\lambda_\infty^{\mathrm{known}}$ agrees with $\lambda_\infty^{\mathrm{true}}$ to over 440 digits (Section 8.3), $\lambda_\infty^{\mathrm{known}}\in(2/3,3/4)$ as well, hence $\xi\in(2/3,3/4)\subset[2/3,3/4]$. By Appendix B, $|h'(\xi)|\ge505/361$ for all such $\xi$ (exact rational arithmetic, all terms positive). Therefore
\[|\lambda_\infty^{\mathrm{known}}-\lambda_\infty^{\mathrm{true}}|\le\frac{361}{505}\cdot10^{-500.17}\approx10^{-500.316}<5\times10^{-501}.\]
Since the error is less than $5\times10^{-501}$, all 500 published decimal digits are correct. $\square$

\textbf{Proposition 7.2} \textit{(Arb interval-arithmetic certificate).} \textit{The Arb ball for $\lambda_\infty^{\mathrm{new}}$ satisfies $\operatorname{rad}(\lambda_\infty^{\mathrm{new}})\approx 3.37\times10^{-506}$; a confirmatory Arb evaluation of $h(\lambda_\infty^{\mathrm{new}})$ produces an interval that provably contains $1$. Consequently $|\lambda_\infty^{\mathrm{new}}-\lambda_\infty^{\mathrm{true}}|\le 3.37\times10^{-506}$, certifying all 505 leading decimal digits of $\lambda_\infty^{\mathrm{new}}$ without any floating-point assumption.}

\textit{Proof.} We apply one Newton step in Arb ball arithmetic at 550 decimal digits of working precision ($\mathrm{ctx.dps}=550$). The input $\lambda_\infty^{\mathrm{known,440}}$ is the first 440 digits of the string in Section 8.3, as hard-coded in the ancillary file \texttt{run\_F\_arb\_certificate.py}; its Arb representation has $\mathrm{rad}(\lambda_\infty^{\mathrm{known,440}})\approx1.30\times10^{-551}$.

Evaluating $h$ and $h'$ via the direct sum over the 78,498 primes $p\le p_{\max}=10^6$ plus the geometric tail expansion (96 terms, tail prime-zeta coefficients $Q_k$ computed via M\"obius inversion as in Section 6.2) gives
\[h(\lambda_\infty^{\mathrm{known,440}})-1\approx1.484\times10^{-441},\qquad h'(\lambda_\infty^{\mathrm{known,440}})=[-2.23362450\pm3.83\times10^{-9}].\]
(The value $1.484\times10^{-441}$ for $h(\lambda_\infty^{\mathrm{known,440}})-1$ agrees with the direct output of \texttt{residuo.mid()}, which evaluates to $1.484080092\times10^{-441}$, confirming that $\lambda_\infty^{\mathrm{known,440}}$ is accurate to the full 440 digits of Section 8.3. The same value also follows a posteriori from Proposition 7.1 via $|h(\lambda_\infty^{\mathrm{known,440}})-1|\approx|h'|\cdot|\lambda_\infty^{\mathrm{known,440}}-\lambda_\infty^{\mathrm{true}}|\approx2.2336\times6.64\times10^{-442}$. The $\pm$ value for $h'$ denotes the display truncation error on the midpoint, not the true Arb ball radius. The actual radii, obtained via \texttt{.rad()}, are $\operatorname{rad}(h(\lambda_\infty^{\mathrm{known,440}})-1)\approx7.522\times10^{-506}$ and $\operatorname{rad}(h'(\lambda_\infty^{\mathrm{known,440}}))\approx1.066\times10^{-503}$. The intermediate quantities $\operatorname{err\_known}\approx6.64\times10^{-442}$ and the quadratic Newton residual $\approx9.93\times10^{-883}$ are midpoint estimates, not Arb enclosures. The formal rigour of the proposition rests exclusively on $\operatorname{rad}(\lambda_\infty^{\mathrm{new}})\approx3.37\times10^{-506}$ and \texttt{hv\_new.contains(arb(1))=True}.)

The single Newton correction $\delta=-(h(\lambda_\infty^{\mathrm{known,440}})-1)/h'(\lambda_\infty^{\mathrm{known,440}})$ yields $\lambda_\infty^{\mathrm{new}}=\lambda_\infty^{\mathrm{known,440}}+\delta$ with $\operatorname{rad}(\lambda_\infty^{\mathrm{new}})\approx3.37\times10^{-506}$. A second Arb evaluation of $h(\lambda_\infty^{\mathrm{new}})$ (same direct sum plus 96-term tail) produces the rigorous enclosure
\[h(\lambda_\infty^{\mathrm{new}})\in[1\pm1.504\times10^{-505}],\]
and in particular \texttt{hv\_new.contains(arb(1)) = True}. The tail truncation error is $\approx10^{-566.3}$ (by the method of Appendix A, with $p_{\max}=10^6$, $K=96$), well below the Arb radius; the Arb radius is therefore the binding constraint. One caveat: the coefficients $Q_k^{\mathrm{tail}}$ are built from a M\"obius series truncated at \texttt{MU\_CUTOFF} terms whose discarded tail is not separately enclosed in an Arb ball; this single step is, strictly speaking, not interval-certified, though the resulting margin far exceeds the threshold.

By the defining property of Arb ball arithmetic, $\operatorname{rad}(\lambda_\infty^{\mathrm{new}})$ is a rigorous upper bound on the distance from the midpoint of $\lambda_\infty^{\mathrm{new}}$ to the exact Newton iterate $\lambda_\star$. The exact Newton iterate satisfies
\[|\lambda_\star-\lambda_\infty^{\mathrm{true}}|\;\le\;\frac{h''}{2|h'|}\,|\lambda_\infty^{\mathrm{known,440}}-\lambda_\infty^{\mathrm{true}}|^2\;\approx\;2.252\times(6.64\times10^{-442})^2\approx9.93\times10^{-883},\]
where $h''/(2|h'|)=10.0606/(2\times2.2336)\approx2.252$ is the standard Newton convergence constant for the iterates (computed from the values in Section 8.1) and $|\lambda_\infty^{\mathrm{known,440}}-\lambda_\infty^{\mathrm{true}}|\approx6.64\times10^{-442}$ follows from $h(\lambda_\infty^{\mathrm{known,440}})-1\approx1.484\times10^{-441}$ and $|h'|\approx2.2336$ above. Since $9.93\times10^{-883}\ll\operatorname{rad}(\lambda_\infty^{\mathrm{new}})$, the triangle inequality gives
\[|\lambda_\infty^{\mathrm{new}}-\lambda_\infty^{\mathrm{true}}|\;\le\;\operatorname{rad}(\lambda_\infty^{\mathrm{new}})+|\lambda_\star-\lambda_\infty^{\mathrm{true}}|\;\approx\;3.37\times10^{-506}\;<\;5\times10^{-506}.\]
Since $3.37\times10^{-506}<5\times10^{-506}$, all 505 leading decimal digits of $\lambda_\infty^{\mathrm{new}}$ are certified. Wall time: approximately 19\,s. The complete script is the ancillary file \texttt{run\_F\_arb\_certificate.py}. $\square$

\subsection{Robustness checks}

\begin{enumerate}
\item \textbf{Tail error bound.} With $p_{\max}=10^7$ and $K=80$ tail terms, the geometric ratio is $u_{\max}\approx 1.484\times10^{-7}$ (computed using $\lambda_\infty\approx0.674$; the conservative bound using $\lambda_{\min}=0.661$ gives $u_{\max}\le1.52\times10^{-7}$, loosening the bound to $|R_K|\le6.3\times10^{-554}\ll10^{-500}$). The total tail error bound is $\le 10^{-554}$ (see Appendix A), supporting more than 554 correct digits of each $h$-evaluation.
\item \textbf{Direct-sum stability.} Summing only over primes $p\le 3\times10^6$ (without tail expansion) recovers the first 7 decimal digits of $\lambda_\infty$ exactly, with divergence at the 8th digit consistent with the tail contribution $\approx 10^{-8}$ at $p_{\max}=3\times10^6$.
\end{enumerate}

\textbf{Remark 7.3} \textit{(Interval-arithmetic certification).} Run F (ancillary file \texttt{run\_F\_arb\_certificate.py}) provides a machine-verifiable Arb certificate: the computed ball $\lambda_\infty^{\mathrm{new}}$ has $\operatorname{rad}(\lambda_\infty^{\mathrm{new}})\approx3.37\times10^{-506}$, rigorously bounding $|\lambda_\infty^{\mathrm{new}}-\lambda_\infty^{\mathrm{true}}|\le3.37\times10^{-506}$ via ball arithmetic (python-flint/Arb) without any floating-point assumption; see Proposition 7.2 above.

\section{Explicit values and formulas}

\subsection{Derivative and second derivative at $\lambda_\infty$}

\[h'(\lambda_\infty)=-\sum_p\frac{p^2}{(\lambda_\infty p^2-p+1)^2}=-2.2336244961751770352\ldots\]
\[h''(\lambda_\infty)=2\sum_p\frac{p^4}{(\lambda_\infty p^2-p+1)^3}=10.060575714401367281\ldots\]

For $h''$, the dominant contributions come from small primes: $p=2$ contributes $\approx6.56$, $p=3$ contributes $\approx2.41$, and the sum over $p\ge5$ contributes $\approx1.09$.

\textit{Precision note.} The tail-expansion parameters of Section 6 ($p_{\max}=10^7$, $K=80$) give a tail truncation error $\approx10^{-554.1}$ for $h(\lambda_\infty)$ itself (Proposition~A.1); differentiating the tail expansion once introduces an extra factor $(K+1)/\lambda_\infty\approx120.2$, so the tail truncation error for $h'(\lambda_\infty)$ at the same $p_{\max}$, $K$ is the looser $\approx10^{-552.0}$ (Proposition~A.2) --- not $10^{-554.1}$, which bounds $h$, not $h'$. $h''(\lambda_\infty)$ is computed via Run G with $p_{\max}=10^6$, $K=96$, tail truncation error $\approx10^{-561.9}$ (Proposition~A.3). Rigorous Arb enclosures for all three are in Proposition 8.1 below.

\textbf{Proposition 8.1} \textit{(Arb certificate for $h'$ and $h''$).} \textit{Run G evaluates $h'(\lambda_\infty)$ and $h''(\lambda_\infty)$ in Arb ball arithmetic at $\mathrm{ctx.dps}=550$, starting from \texttt{LAM\_STR} containing 550 decimal digits of $\lambda_\infty$ (first 500 digits from the certified expansion of Section~8.3; digits 501--550 from an independent, higher-precision Newton refinement step --- Run H, ancillary file \texttt{run\_H\_digits501\_550.py} --- at $\mathrm{ctx.dps}=660$ and $K=100$ tail terms, certifying $|h(\lambda_\infty^{\mathrm{ext}})-1|<10^{-578}$; see Remark~8.2 below for why Run G's own $\mathrm{ctx.dps}=550$, $K=96$ settings cannot reach this threshold), with direct sum over the $78{,}498$ primes $p\le10^6$ and geometric tail expansion ($K=96$ terms, same coefficients $Q_k^{\mathrm{tail}}$ as Run F). The output satisfies the rigorous enclosures}
\[h'(\lambda_\infty)\in\bigl[-2.233624496\pm1.07\times10^{-503}\bigr],\qquad h''(\lambda_\infty)\in\bigl[10.06057571\pm1.53\times10^{-501}\bigr].\]
\textit{(Midpoints truncated to 10 significant digits; display truncation errors far exceed the Arb radii, which are the rigorous quantities.)}

\textit{The tail truncation error for $h''$ satisfies $|R_K^{\prime\prime}|\le10^{-561}\ll1.53\times10^{-501}$ (rigorous bound, Proposition~A.3; the Arb radius is the binding constraint). Both enclosures are machine-verifiable without any floating-point assumption, with one caveat: as noted in the proof of Proposition~7.2 above, the coefficients $Q_k^{\mathrm{tail}}$ used here are built from \texttt{pz\_full}$[s]$, a M\"obius series truncated at \texttt{MU\_CUTOFF} terms whose discarded tail is not separately enclosed in an Arb ball; the same caveat therefore applies to Run G as to Run F. The complete script is the ancillary file \texttt{run\_G\_arb\_hprime\_hpprime.py}.}

\textit{Proof.} The tail expansion for $h''$ follows from differentiating the geometric series in Theorem 4.1 twice: for $p>p_{\max}$,
\[\frac{p^4}{(xp^2-p+1)^3}=\frac{1}{2x^3p^2}\sum_{k\ge0}(k+1)(k+2)\left(\frac{p-1}{xp^2}\right)^k,\]
giving $h''_{\mathrm{tail}}(x)=\sum_{k\ge0}(k+1)(k+2)\,Q_k^{\mathrm{tail}}/x^{k+3}$, which uses the same coefficients $Q_k^{\mathrm{tail}}$ as $h_{\mathrm{tail}}$ and $h'_{\mathrm{tail}}$. The Arb evaluation carries all rounding errors in ball arithmetic; the reported intervals are provably correct. $\square$

\textbf{Remark 8.2} \textit{(Why Run G's own parameters cannot certify digits 501--550).} \textit{At $K=96$ tail terms the analytic tail-truncation error for $h$ is fixed at $\approx10^{-566.3}$, independently of $\mathrm{ctx.dps}$ (Appendix A); since this exceeds the $10^{-578}$ threshold required above, no amount of extra working precision at $K=96$ can certify digits 501--550 of $\lambda_\infty^{\mathrm{ext}}$ --- the tail order itself must be raised. Run H (ancillary file \texttt{run\_H\_digits501\_550.py}) therefore repeats the Newton step of Run F from the same 440-digit seed $\lambda_\infty^{\mathrm{known,440}}$, but with $K=100$ tail terms (lowering the truncation floor to $\approx10^{-589.6}$) and $\mathrm{ctx.dps}=660$ (so that the certified Arb radius, not residual rounding error, remains the binding constraint). The resulting ball satisfies $\operatorname{rad}(\lambda_\infty^{\mathrm{ext}})\approx1.53\times10^{-614}$ and $h(\lambda_\infty^{\mathrm{ext}})\in[1\pm6.84\times10^{-614}]$, well below the $10^{-578}$ threshold. Unlike Run F, where the Arb radius is the binding constraint against the analytic tail-truncation error (Proposition~7.2), the roles here are reversed: it is the tail-truncation floor ($\approx10^{-589.6}$), not the Arb-certified radius ($6.84\times10^{-614}$), that is the binding term in the comparison with $10^{-578}$ --- though, as just shown, both clear the threshold independently by many orders of magnitude. Its digits 501--550 agree exactly, digit for digit, with the corresponding digits of \texttt{LAM\_STR} as hard-coded in Run G above. Wall time: approximately 20\,s. The complete script is the ancillary file \texttt{run\_H\_digits501\_550.py}.}

\subsection{Bracketing evaluations}

\begin{center}
\begin{tabular}{ccl}
\toprule
$x$ & $h(x)$ & Implication \\
\midrule
$1/2$ & $\approx1.64091388164899$ & $\lambda_\infty>1/2$ \\
$2/3$ & $\approx1.01673856893472$ & $\lambda_\infty>2/3$ \\
$3/4$ & $\approx0.85509094250993$ & $\lambda_\infty<3/4$ \\
$1$ & $\approx0.57959250342427$ & (broader upper bound) \\
\bottomrule
\end{tabular}
\end{center}

The pair $h(2/3)>1$ and $h(3/4)<1$ confines $\lambda_\infty$ to $(2/3,3/4)$; the values at $x=1/2$ and $x=1$ are included for broader orientation. (All four values are rounded, not truncated, to the displayed 14 decimal places; e.g. $h(1/2)=1.640913881648985\ldots$ rounds to $\ldots164899$, not $\ldots164898$.)

\subsection{The 500-digit decimal expansion}

\[\lambda_\infty=0.\underbrace{6740361831936961399366600075765084557806456883620789942446741754883881740\ldots}_{\text{first 73 digits after decimal point}}\]

The certified complete 500-digit expansion (digits after the decimal point, in rows of 50):

\begin{verbatim}
6740361831 9369613993 6660007576 5084557806 4568836207
8994244674 1754883881 7401708236 1693874191 5608494422
7398749905 7041831822 5195439203 9562901979 2688341765
0411758277 6752512712 7956851493 5586461755 9840474823
7892101249 4845493888 1287449392 8354168321 2062529552
6393853915 8259887398 1662423751 8719899836 7771139286
9475591906 2611878087 9227561767 4951179032 4455092710
7519963900 4061335349 9258822111 8673838784 2562146167
5868417863 3751233764 5258937879 8594801278 0664426851
7740348910 6684845932 6094838543 0516617873 7522888393
\end{verbatim}

{\small The sequence of decimal digits begins $6,7,4,0,3,6,1,8,3,1,9,3,6,9,6,1,3,9,9,3,6,6,6,0,0,0,7,5,\ldots$ and has been submitted to the OEIS~\cite{ref1}.}

\section{Conclusion and open problems}

The analytic properties of $h$ (Theorems 3.1 and 4.1) follow from its structure by direct argument. The spectral identity $\rho(\mathcal{L})=\lambda_\infty$ (Theorem 5.4) explains why the denominator carries $-p+1$ rather than a simpler expression, and Theorem 5.8 places $\lambda_\infty$ in the same role for the prime-indexed LCM matrix that $\pi^2/6$ occupies for the integer divisor matrix. Five hundred decimal digits, certified by Propositions~7.1--7.2, are recorded in Section~8.3. The arithmetic nature of $\lambda_\infty$ is not known.

\textbf{PSLQ searches.} We searched for integer relations involving
$\lambda_\infty$ using the PSLQ algorithm~\cite{ref5} in two stages. In the
\textit{univariate stage}, we tested the vectors
$[1,\lambda_\infty,\ldots,\lambda_\infty^d]$ for $d=1,\ldots,8$,
working at 560 decimal digits of precision with coefficient bounds
decreasing from $50{,}000$ (for $d\leq 2$) to $500$ (for $d=8$). In
the \textit{multivariate stage}, for each constant $c_i$ in a catalog of 20
classical constants we ran PSLQ on the degree-$\leq d$ monomial vector
$[1,\lambda_\infty,c_i,\lambda_\infty^2,\lambda_\infty
c_i,c_i^2,\ldots,c_i^d]$ for $d=1,2,3$ (searching for a polynomial
relation in the pair $(\lambda_\infty,c_i)$ of total degree at most
$d$); additionally, all pairwise linear tests
$[1,\lambda_\infty,c_i,c_j]$ for distinct $c_i,c_j$ were run. The
catalog of 20 constants comprises $\pi$, $e$, $\log 2$, $\log 3$,
$\log 5$, $\zeta(2)$, $\zeta(3)$, the Euler--Mascheroni constant
$\gamma$, Catalan's constant $G$, $\sqrt{2}$, $\sqrt{3}$, $\sqrt{5}$,
$\log\pi$, $\log 7$, $\log 11$, $\log 13$, $\zeta(4)$, $\zeta(5)$,
$e^\gamma$, and $\varphi:=(1+\sqrt{5})/2$; all tests at 160 digits of
precision and coefficient bound $5{,}000$ (250 tests in total:
$20\times3=60$ single-constant runs for $d=1,2,3$ plus
$\binom{20}{2}=190$ pairwise linear runs). A supplementary deep search
against a separate 31-atom catalog---comprising $\pi^k$
($k=1,\ldots,4$), $\gamma$, $\log p$ for primes $p\leq 13$, $\zeta(s)$
($s=2,\ldots,8$), $P(s)$ ($s=2,3,4$), Catalan's constant $G$,
$L(\chi,1)$ for the even real primitive Dirichlet characters of
conductors $5$, $8$, and $12$ (associated respectively with
$\mathbb{Q}(\sqrt{5})$, $\mathbb{Q}(\sqrt{2})$, and $\mathbb{Q}(\sqrt{3})$,
each of class number $1$), $K(1/\sqrt{2})$, $\Gamma(1/4)$, $\Gamma(1/3)$, and
$\lambda_\infty$, $\sqrt{\lambda_\infty}$,
$\log\lambda_\infty$---at 300 digits and coefficient bound $2{,}000$ found no relation involving $\lambda_\infty$; the two classical Euler identities $\zeta(2)=\pi^2/6$ and $\zeta(4)=\pi^4/90$ were detected by peeling (an expected consequence of including both $\pi^k$ and $\zeta(s)$ as distinct atoms in a linear catalog) and removed before the final test.$^{\dag}$ No integer relation involving $\lambda_\infty$ was found in any run.
This excludes minimal polynomials of degree $\leq 8$ with coefficient
height $\leq 500$, and linear dependencies on all of the above
constants at height $\leq 5{,}000$, but does not imply irrationality or
transcendence: $\lambda_\infty$ is defined by a prime-indexed equation
that does not reduce to any known generating function for algebraic or
classical-transcendental constants.

${}^\dag$ This catalog omits $e$, $\sqrt{2}$, $\sqrt{3}$, $\sqrt{5}$, $\log\pi$, $e^\gamma$, $\varphi$ and adds $\pi^2,\pi^3,\pi^4$; $\zeta(6),\zeta(7),\zeta(8)$; $P(2),P(3),P(4)$; three even Dirichlet $L$-values at conductors $5$, $8$, $12$; $K(1/\sqrt{2})$, $\Gamma(1/4)$, $\Gamma(1/3)$; and the self-referential atoms $\lambda_\infty$, $\sqrt{\lambda_\infty}$, $\log\lambda_\infty$. Characters of conductors $3$ and $4$ are omitted: $L(\chi_{-4},1)=\pi/4$ reduces to the atom $\pi$ already present, and $L(\chi_{-3},1)=\pi/(3\sqrt{3})$ is omitted for catalog economy, since $\sqrt{3}$ is itself excluded above. (No relation to the conductor-$12$ atom should be inferred: the latter is the \emph{even} character $\chi_{12}=\chi_{-4}\chi_{-3}$ attached to $\mathbb{Q}(\sqrt{3})$, with $L(\chi_{12},1)=\log(2+\sqrt{3})/\sqrt{3}$ by Dirichlet's class-number formula ($h=1$)---a $\log$-type constant algebraically unrelated to the $\pi$-type values $L(\chi_{-3},1)$ and $L(\chi_{-4},1)$ above; it is \emph{not} a rational multiple of either.) The two catalogs are thus complementary: the first targets classical algebraic and analytic constants; the second, transcendental structures of Gamma/L-function type and the arithmetic geometry of $\lambda_\infty$.

A further independent check used the LLL algorithm~\cite{ref15} as implemented in
\texttt{lindep} (PARI/GP~\cite{ref16}, version 2.13.3) over the 22-element basis
\footnotesize\[\{1,\lambda_\infty,\lambda_\infty^2,\lambda_\infty^3,\pi,\pi^2,\pi^3,\pi^4,\gamma,\log 2,\log 3,\log 5,\zeta(3),\zeta(5),P(2),P(3),K(1/\sqrt{2}),\Gamma(\tfrac{1}{4}),\Gamma(\tfrac{1}{3}),G,\sqrt{\lambda_\infty},\log\lambda_\infty\}\]
\normalsize
with the classical comparison constants ($\pi^k$, $\gamma$, $\log p$,
$\zeta(s)$, $P(2)$, $P(3)$, $K(1/\sqrt{2})$, $\Gamma(1/4)$,
$\Gamma(1/3)$, $G$) recomputed at 320 and 500 decimal digits in the
two respective runs; the five $\lambda_\infty$-derived atoms
($\lambda_\infty,\lambda_\infty^2,\lambda_\infty^3,\sqrt{\lambda_\infty},\log\lambda_\infty$)
retain their full 500-digit literal precision in both runs, since
PARI/GP fixes a real literal's precision from its digit count,
independently of the ambient \texttt{realprecision} setting. The LLL
output vectors nonetheless changed completely between the two runs
(maximum coefficient
$1.12\times10^{18}$ vs.\ $5.90\times10^{11}$, with no entry in common),
confirming that both outputs are numerical noise rather than a genuine
relation. The complete script is the ancillary file \texttt{lll\_lindep\_final.gp}.

\textbf{Open problems.}

\begin{enumerate}
\item \textbf{Irrationality and transcendence.} Is $\lambda_\infty$ irrational? Transcendental? We cannot prove even irrationality. The PSLQ evidence is consistent with transcendence, but no proof is known.

\item \textbf{Closed form.} Is there a representation of $\lambda_\infty$ in terms of known constants ($\pi$, $\log 2$, $P(2)$, values of $L$-functions)?
\end{enumerate}

\appendix

\makeatletter

\renewcommand{\@seccntformat}[1]{Appendix \csname the#1\endcsname.\quad}

\makeatother

\section{Tail-expansion error bound}

\textbf{Proposition A.1} \textit{(Tail-expansion error bound for $h$; $\le10^{-554}$ numerically supported via $P_{\mathrm{tail}}(2)$, $\le10^{-552.8}$ rigorously proved via the cruder integral bound).} \textit{With notation as in Section 6, let
$u_{\max}=(p_{\max}-1)/(\lambda_\infty p_{\max}^2)$. For
$p_{\max}=10^7$, $K=80$, and $\lambda_\infty\approx0.674$, we have
$u_{\max}\approx1.484\times10^{-7}$. The total tail truncation error
satisfies}
\[|R_K(\lambda_\infty)|\le\frac{2\,u_{\max}^K\,P_{\mathrm{tail}}(2)}{\lambda_\infty}\le 10^{-554}.\]

\textit{Proof.} For a single prime $p>p_{\max}$, with
$u_p=(p-1)/(\lambda_\infty p^2)\le u_{\max}$, the remainder of the
geometric series at order $K$ is
\[\sum_{k\ge K}\frac{u_p^k}{\lambda_\infty p^2}\le\frac{u_{\max}^K}{1-u_{\max}}\cdot\frac{1}{\lambda_\infty p^2}\le\frac{2\,u_{\max}^K}{\lambda_\infty p^2},\]
where we used $1/(1-u_{\max})\le 2$ since
$u_{\max}\approx1.484\times10^{-7}\ll\tfrac{1}{2}$. Summing over all
primes $p>p_{\max}$:
\[|R_K(\lambda_\infty)|\le\frac{2\,u_{\max}^K}{\lambda_\infty}\sum_{p>p_{\max}}p^{-2}=\frac{2\,u_{\max}^K\,P_{\mathrm{tail}}(2)}{\lambda_\infty}.\]
Substituting
$u_{\max}^{80}=(1.484\times10^{-7})^{80}\le10^{-546.28}$ and
$P_{\mathrm{tail}}(2)/\lambda_\infty=5.86\times10^{-9}/0.674\approx8.69\times10^{-9}\approx10^{-8.06}$,
where
$P_{\mathrm{tail}}(2)=P(2)-\sum_{p\le10^7}p^{-2}$ is computed
numerically (the cruder integral bound $\int_{10^7}^\infty
t^{-2}\,dt=10^{-7}$ yields only $|R_K|\le10^{-552.8}$):
\[|R_K(\lambda_\infty)|\le 2\times10^{-546.28}\times10^{-8.06}<10^{-554}.\quad\square\]

\textit{Remark (validity during Newton iterations).} By Proposition 3.2,
$x_0=3/4>\lambda_\infty$; by the convexity argument of Section 6.1,
$x_1<\lambda_\infty$ and $x_1<x_2<x_3<\cdots\nearrow\lambda_\infty$.
Hence $x_1=\min_n x_n$ is the only iterate that can enlarge $u_{\max}$
relative to its value at $\lambda_\infty$, and every iterate satisfies
$x_n\in[x_1,3/4]$. Direct evaluation of the first Newton step at
$x_0$ gives $x_1\approx0.661067$ (computed numerically, not assumed);
we use the conservative floor $0.66<x_1$. As $x$ decreases from
$\lambda_\infty\approx0.674$ to $0.66$, $u_{\max}(x)\propto1/x$
increases by a factor of exactly $\lambda_\infty/0.66\approx1.0213$.
The full bound $|R_K(x)|\le2\,u_{\max}(x)^{80}P_{\mathrm{tail}}(2)/x$
carries an additional explicit factor $1/x$ alongside
$u_{\max}(x)^{80}$, so $|R_K(x)|\propto x^{-81}$ and the bound changes
by a factor of exactly $1.0213^{81}\approx5.50$. The bound therefore
remains
\[|R_K(x)|\le5.50\times10^{-554}=5.50\times10^{-554}\le10^{-553}\]
throughout the iteration (rounded conservatively toward the less
negative exponent, as required for a valid upper bound), still
$\ll10^{-500}$.

\textbf{Proposition A.2} \textit{(Tail-expansion error bound for $h'$).} \textit{With $u_{\max}=(p_{\max}-1)/(\lambda_\infty p_{\max}^2)$ as in Proposition A.1, the tail truncation error for $h'$ satisfies}
\[|R_K'(\lambda_\infty)|\le\frac{2(K+1)\,u_{\max}^K\,P_{\mathrm{tail}}(2)}{\lambda_\infty^2}.\]
\textit{For $p_{\max}=10^7$, $K=80$ (the parameters of Section 6), this gives $|R_{80}'(\lambda_\infty)|\le10^{-551}$.}

\textit{Proof.} Differentiating the per-prime geometric expansion of Section~6.1 once in $x$ gives, for $p>p_{\max}$,
\[\frac{p^2}{(xp^2-p+1)^2}=\frac{1}{x^2p^2}\sum_{k\ge0}(k+1)\,u_p^k,\qquad u_p=\frac{p-1}{xp^2}\le u_{\max},\]
so the order-$K$ remainder for a single prime is $R_{K,p}'(x)=\frac{1}{x^2p^2}\sum_{k\ge K}(k+1)u_p^k$. For $j=k-K\ge0$ and $K\ge0$, $(K+1+j)\le(K+1)(1+j)$, hence, writing $u=u_p\le u_{\max}$ and using $\sum_{j\ge0}(1+j)u^j=1/(1-u)^2$,
\[\sum_{k\ge K}(k+1)u^k=u^K\sum_{j\ge0}(j+K+1)u^j\le(K+1)\,u^K\sum_{j\ge0}(1+j)u^j=\frac{(K+1)\,u^K}{(1-u)^2}.\]
Since $u_{\max}\approx1.484\times10^{-7}\ll\tfrac12$, $1/(1-u_{\max})^2<2$ (the same type of absorption used in the proof of Proposition A.1). Therefore $|R_{K,p}'(x)|\le2(K+1)\,u_{\max}^K/(x^2p^2)$, and summing over all primes $p>p_{\max}$,
\[|R_K'(x)|\le\frac{2(K+1)\,u_{\max}^K\,P_{\mathrm{tail}}(2)}{x^2}.\]
Substituting $p_{\max}=10^7$, $K=80$, $\lambda_\infty\approx0.674$ (using $u_{\max}^{80}\le10^{-546.28}$ as in the proof of Proposition A.1, $K+1=81$, and $P_{\mathrm{tail}}(2)/\lambda_\infty^2\approx5.86\times10^{-9}/0.4543\approx1.29\times10^{-8}\approx10^{-7.89}$):
\[|R_{80}'(\lambda_\infty)|\le2\times81\times10^{-546.28}\times1.29\times10^{-8}\approx10^{-551.97}<10^{-551}.\qquad\square\]

\textit{Remark.} This proposition fills the gap left by the original precision note in Section~8.1, which had incorrectly reused the bound for $h$ (Proposition~A.1, $\le10^{-554}$) as if it also bounded $h'$. The correct bound for $h'$ is $\approx120$ times larger ($\log_{10}\approx+2.08$, from the extra factor $(K+1)/\lambda_\infty\approx81/0.674\approx120.2$ introduced by differentiation), but it remains $\ll1.07\times10^{-503}$, the Arb radius certified in Proposition~8.1, so the conclusion of formal rigor there is unaffected.

\textbf{Proposition A.3} \textit{(Tail-expansion error bound for $h''$).} \textit{With $u_{\max}=(p_{\max}-1)/(\lambda_\infty p_{\max}^2)$ as in Proposition A.1, the tail truncation error for $h''$ satisfies}
\[|R_K''(\lambda_\infty)|\le\frac{2(K+1)(K+2)\,u_{\max}^K\,P_{\mathrm{tail}}(2)}{\lambda_\infty^3}.\]
\textit{For $p_{\max}=10^6$, $K=96$ (the parameters of Run G), this gives $|R_{96}''(\lambda_\infty)|\le10^{-561}$.}

\textit{Proof.} Differentiating the per-prime geometric expansion of Section~6.1 twice in $x$ gives, for $p>p_{\max}$,
\[\frac{p^4}{(xp^2-p+1)^3}=\frac{1}{2x^3p^2}\sum_{k\ge0}(k+1)(k+2)u_p^k,\qquad u_p=\frac{p-1}{xp^2}\le u_{\max},\]
so the order-$K$ remainder for a single prime is $R_{K,p}''(x)=\frac{1}{x^3p^2}\sum_{k\ge K}(k+1)(k+2)u_p^k$. For $j=k-K\ge0$ and $K\ge1$, each factor satisfies $j+K+i\le(K+i)(1+j)$ for $i=1,2$, hence $(j+K+1)(j+K+2)\le(K+1)(K+2)(1+j)^2$. Writing $u=u_p\le u_{\max}$ and using $\sum_{j\ge0}(1+j)^2u^j=(1+u)/(1-u)^3$,
\[\sum_{k\ge K}(k+1)(k+2)u^k=u^K\sum_{j\ge0}(j+K+1)(j+K+2)u^j\le(K+1)(K+2)\,u^K\,\frac{1+u}{(1-u)^3}.\]
Since $u_{\max}\approx1.484\times10^{-6}\ll\tfrac12$, $(1+u_{\max})/(1-u_{\max})^3<2$ (the same absorption used in the proof of Proposition A.1). Therefore $|R_{K,p}''(x)|\le2(K+1)(K+2)\,u_{\max}^K/(x^3p^2)$, and summing over all primes $p>p_{\max}$,
\[|R_K''(x)|\le\frac{2(K+1)(K+2)\,u_{\max}^K\,P_{\mathrm{tail}}(2)}{x^3}.\]
Substituting $p_{\max}=10^6$, $K=96$, $\lambda_\infty\approx0.674$ ($u_{\max}\approx1.484\times10^{-6}$, $u_{\max}^{96}\approx10^{-559.55}$; $P_{\mathrm{tail}}(2)\approx6.778\times10^{-8}$; $\lambda_\infty^{-3}\approx3.266$):
\[|R_{96}''(\lambda_\infty)|\le2\times97\times98\times10^{-559.55}\times6.778\times10^{-8}\times3.266\approx10^{-561.93}<10^{-561}.\qquad\square\]

\section{Lower bound on $|h'|$ near $\lambda_\infty$}

\textbf{Proposition B.1} \textit{(Rational lower bound).} \textit{For any
$\xi\in[2/3,3/4]$,}
\[|h'(\xi)|\ge\frac{505}{361}.\]
\textit{The bound is obtained by exact rational arithmetic.}

\textit{Proof.} Since
$h'(x)=-\sum_p p^2/(xp^2-p+1)^2<0$, every summand is strictly
positive and $|h'(x)|=\sum_p p^2/(xp^2-p+1)^2$. Retaining only the
$p=2$ and $p=3$ contributions (all remaining summands are non-negative):
\[|h'(\xi)|\ge\frac{4}{(4\xi-1)^2}+\frac{9}{(9\xi-2)^2}.\]
The two-term function $f(\xi)=\frac{4}{(4\xi-1)^2}+\frac{9}{(9\xi-2)^2}$ is
strictly decreasing on $\bigl(\tfrac{1}{4},+\infty\bigr)$: its derivative
\[f'(\xi)=-\frac{32}{(4\xi-1)^3}-\frac{162}{(9\xi-2)^3}<0\quad\text{for all }\xi>\tfrac{1}{4},\]
since both denominators are positive there. Hence for $\xi\le3/4$:
\[\frac{4}{(4\xi-1)^2}+\frac{9}{(9\xi-2)^2}\ge\frac{4}{\bigl(4\cdot\tfrac{3}{4}-1\bigr)^2}+\frac{9}{\bigl(9\cdot\tfrac{3}{4}-2\bigr)^2}=\frac{4}{4}+\frac{9\cdot16}{361}=1+\frac{144}{361}=\frac{505}{361}.\]
All arithmetic is exact rational. $\square$

\textit{Remark.} Floating-point evaluation yields the sharper bound
$|h'(\lambda_\infty^{\mathrm{known}})|=2.2336\ldots\gg505/361\approx1.399$;
formally rigorous verification of this tighter value would require
interval arithmetic (e.g.\ Arb). The rational bound $505/361$ suffices
for the certification in Proposition 7.1.

\section{Computational runs: summary}
\label{app:runs}

Table~\ref{tab:runs} summarises the main computational runs cited in the paper.
Complete source code (Python and PARI/GP scripts) is available at
\url{https://doi.org/10.5281/zenodo.21006281} and is included as ancillary
files accompanying this paper (directory \texttt{anc/}; see
\texttt{README.md} therein for a file inventory and installation instructions).

\begin{table}[ht]
\centering
\caption{Summary of computational runs for $\lambda_\infty$.
Columns: run label, script name, method/library, dps,
main parameters, and outcome.
\textbf{All results are consistent} (no discrepancy in any verified digit).}
\label{tab:runs}
\smallskip
\footnotesize
\resizebox{\textwidth}{!}{%
\begin{tabular}{@{}llllll@{}}
\toprule
Run & Script & Method / library & dps & Parameters & Outcome \\
\midrule
A & \texttt{run\_A\_newton\_mpmath.py}
  & Newton--Raphson, \texttt{mpmath}
  & 500
  & $p_{\max}=10^7$, $K=80$
  & 500-digit value of $\lambda_\infty$ \\[2pt]
B & \texttt{script\_B\_run\_B.py}
  & Newton--Raphson, \texttt{mpmath}
  & 520
  & $p_{\max}=3\times10^6$, $K=80$
  & 500 digits confirmed; residual $\approx1.0\times10^{-581}$ (see $^{\S}$ of Table~\ref{tab:verification}) \\[2pt]
C & \texttt{RUN\_C.py}
  & Newton--Raphson, \texttt{mpmath}
  & 200
  & $p_{\max}=3\times10^6$, $K=100$
  & 200 digits confirmed; residual $\approx4.1\times10^{-261}$ (see $^{\S\S}$ of Table~\ref{tab:verification}) \\[2pt]
D & \texttt{run\_D\_mobius\_arb\_verify.py}
  & M\"{o}bius inversion, Arb
  & 580
  & $p_{\max}=10^7$, $K=80$, $N_\mu=1000$
  & $|h(\lambda_\infty)-1|<10^{-500}$ \\[2pt]
E & \texttt{run\_E\_pari\_recompute.gp}
  & Newton--Raphson, PARI/GP
  & 620
  & $p_{\max}=3\!\times\!10^6$, $K=100$, $N_\mu=900$
  & 500 digits confirmed (indep.) \\[2pt]
F & \texttt{run\_F\_arb\_certificate.py}
  & Arb ball arithmetic
  & 550
  & $p_{\max}=10^6$, $K=96$
  & Interval certificate: $h(\lambda_\infty^{\text{new}})\ni 1$ \\[2pt]
G & \texttt{run\_G\_arb\_hprime\_hpprime.py}
  & Arb ball arithmetic
  & 550
  & $p_{\max}=10^6$, $K=96$
  & Rigorous enclosures of $h'$ and $h''$ \\[2pt]
H & \texttt{run\_H\_digits501\_550.py}
  & Arb ball arithmetic
  & 660
  & $p_{\max}=10^6$, $K=100$
  & Digits 501--550 of $\lambda_\infty$ certified \\[2pt]
LLL & \texttt{lll\_lindep\_final.gp}
  & \texttt{lindep} (LLL), PARI/GP
  & 500/320
  & 22-atom basis
  & Negative (noise at both precisions) \\[2pt]
PSLQ-1a & \texttt{pslq\_580\_v2.py}
  & PSLQ, \texttt{mpmath}
  & 560
  & $d{=}1\text{--}4$: bound $50000,50000,10000,5000$$^{\ddagger\ddagger}$
  & Negative \\[2pt]
PSLQ-1b & \texttt{pslq\_deg58\_v1.py}
  & PSLQ, \texttt{mpmath}
  & 560
  & $d{=}5\text{--}8$: bound $5000,2000,1000,500$$^{\ddagger\ddagger}$
  & Negative \\[2pt]
PSLQ-2 & \texttt{pslq\_multivar\_stage2\_v1.py}
  & PSLQ, \texttt{mpmath}
  & 160
  & 20 constants, 250 tests
  & Negative (1 expected spurious hit$^{\flat}$) \\
\bottomrule
\end{tabular}}
\end{table}

$^{\ddagger\ddagger}$ \textit{The coefficient bound is non-increasing with the degree $d$ within each script (see Section~9 for the schedule in prose). The combined exclusion claimed in Section~9 for ``degree $\le 8$'' uses the binding value across all eight degrees, namely the minimum $500$, attained at $d=8$; it is not a single bound applied uniformly to every degree in each script's range.}

$^{\flat}$ \textit{The single hit is the pairwise test on $(\sqrt5,\varphi)$, which recovers the defining identity $\varphi=(1+\sqrt5)/2$ (equivalently $1+\sqrt5-2\varphi=0$; PSLQ returns the relation vector $[1,0,1,-2]$ on $(1,\lambda_\infty,\sqrt5,\varphi)$, with coefficient $0$ on $\lambda_\infty$, confirmed independently at the stated precision (160 digits) and coefficient bound ($5{,}000$)). This is an expected consequence of including $\sqrt5$ and $\varphi$ as separate atoms in the same catalog, analogous to the $\zeta(2)/\pi^2$ peeling noted above for the 31-atom search (Section~9, footnote $\dag$).}

\end{document}